\pdfoutput=1
\documentclass[10pt]{article}
\usepackage{graphicx}
\usepackage{amssymb}
\usepackage{epstopdf}
\usepackage{amsmath}
\usepackage{amsfonts}

\usepackage{mathrsfs}

\def\phi{{\varphi}}

\def\P{{\mathcal P}}

\DeclareSymbolFont{AMSb}{U}{msb}{m}{n}
\DeclareMathSymbol{\N}{\mathbin}{AMSb}{"4E}
\DeclareMathSymbol{\Z}{\mathbin}{AMSb}{"5A}
\DeclareMathSymbol{\R}{\mathbin}{AMSb}{"52}
\DeclareMathSymbol{\Q}{\mathbin}{AMSb}{"51}
\DeclareMathSymbol{\I}{\mathbin}{AMSb}{"49}
\DeclareMathSymbol{\C}{\mathbin}{AMSb}{"43}

\def\be{\begin{equation}}
\def\ee{\end{equation}}
\def\ber{\begin{eqnarray}}
\def\eer{\end{eqnarray}}

\def\e{\epsilon}

\def\beq{\begin{equation}}
\def\eeq{\end{equation}}

\begin{document}

\addtolength{\textheight}{0 cm} \addtolength{\hoffset}{0 cm}
\addtolength{\textwidth}{0 cm} \addtolength{\voffset}{0 cm}

\newcommand{\ZZ}{\mathbb{Z}}
\newcommand{\Rm}{\mathbb{R}}
\newcommand{\RR}{\mathbb{R}}
\newcommand{\NN}{\mathbb{N}}
\newcommand{\sU}{\mathcal{U}}
\newcommand{\sF}{\mathcal{F}}
\newcommand{\sM}{\mathcal{M}}
\newcommand{\sS}{\mathcal{S}}
\newcommand{\mL}{\mathcal{L}}
\newcommand{\ac}{\hbox{\small ac}}
\newcommand{\mC}{\ensuremath{\mathcal{C}}}
\newcommand{\mU}{\ensuremath{\mathcal{U}}}
\newcommand{\mT}{\ensuremath{\mathcal{T}}}
\newcommand{\mS}{\ensuremath{\mathcal{S}}}
\newcommand{\mF}{\ensuremath{\mathcal{F}}}
\newcommand{\Nm}{\ensuremath{\mathbb{N}}}
\newcommand{\Zm}{\ensuremath{\mathbb{Z}}}
\newcommand{\Hm}{\ensuremath{\mathbb{H}}}
\newcommand{\mM}{\ensuremath{\mathcal{M}}}
\newcommand{\mK}{\ensuremath{\mathcal{K}}}
\newcommand{\mD}{\ensuremath{\mathcal{D}}}
\newcommand{\mA}{\ensuremath{\mathcal{A}}}
\newcommand{\mO}{\ensuremath{\mathcal{O}}}
\newcommand{\mI}{\ensuremath{\mathcal{I}}}
\newcommand{\mB}{\ensuremath{\mathcal{B}}}
\newcommand{\Tm}{\ensuremath{\mathbb{T}}}
\newcommand{\mE}{\ensuremath{\mathcal{E}}}
\newcommand{\vs}{\vspace{.5cm}}

\def\proof {\noindent{\sc{Proof. }}}
\def\qed {\mbox{}\hfill {\small \fbox{}} \\}  
\def\lto{\longrightarrow}
\def\lmto{\longmapsto}
\def\eq{\Longleftrightarrow}
\def\leq{\leqslant}
\def\geq{\geqslant}

\newtheorem{lem}{Lemma}
\newtheorem{thm}{Theorem}
\newtheorem{conj}[lem]{Conjecture}
\newtheorem{ques}[lem]{Question}
\newtheorem{cor}[lem]{Corollary}
\newtheorem{prop}[lem]{Proposition}
\newtheorem{defn}[lem]{Definition}
\newtheorem{note}[lem]{Note}
\newtheorem{rmk}{Remark}
\def\proof {\noindent{\sc{Proof. }}}
\def\qed {\mbox{}\hfill {\small \fbox{}} \\}  

\newcommand{\grad}{\operatorname{grad}}
\newcommand{\Leg}{\mathcal{L}}

\newcommand{\lbstoc}[0]{\underline{B}}
\newcommand{\ubstoc}[0]{\overline{B}^{\text{stoc}}}
\newcommand{\deriv}[2]{\ensuremath{\frac{d{#1}}{d{#2}}}}
\newcommand{\Id}[1]{\ensuremath{\boldsymbol{1}[#1]}}
\newcommand{\abs}[1]{\ensuremath{\left\lvert#1\right\rvert}}
\newcommand{\norm}[1]{\ensuremath{\left\lVert#1\right\rVert}}
\renewcommand{\P}[1]{\ensuremath{\mathbb{P}(#1)}}
\newcommand{\bracket}[1]{\ensuremath{\left(#1\right)}}
\newcommand{\expect}[2][]{\ensuremath{\mathbb{E}_{#1}\left[#2\right]}}
\newcommand{\expcond}[2]{\ensuremath{\mathbb{E}\left[#1\middle\rvert #2\right]}}
\newcommand{\ball}[2]{\ensuremath{B(#1,#2)}}
\newcommand{\proj}[1]{\ensuremath{\text{proj}_{#1}}}
\newcommand{\pderiv}[2]{\ensuremath{\frac{\partial{#1}}{\partial{#2}}}}
\newcommand{\Cmik}[0]{\ensuremath{C}} 
\newcommand{\lip}[1]{\text{Lip}(#1)} 
\newcommand{\dbound}[0]{C^\infty_\text{db}} 
\newcommand{\convex}[0]{c} 
\newcommand{\flip}[1]{\hat{#1}}

\def\pR{\R\cup\{+\infty\}}

\newcommand{\ntx}{\textnormal}
\newcommand\E{{\mathbb E}}

\def\F{{\mathcal {F}}}
\def\e{{\mathcal {E}}}
\def\L{{\mathcal {L}}}
\def\X{{\mathcal {X}}}

\def\msE{\mathscr{E}}
\def\msM{\mathscr{M}}
\def\msH{\mathscr{H}}
\def\msF{\mathscr{F}}
\def\msL{\mathcal{L}}
\def \mbL{\mathbb{L}}
\def\vr{\vert}
\def\Vr{\Vert}

\def\F{{\cal F}}
\def\G{{\cal G}}
\def\L{{\cal L}}
\def\Rm{{\cal R}}
\def\div{\mbox{div} }
\def\t{\tilde}
\def\t{\tilde}

\title{Dynamic and Stochastic Propagation of Brenier's Optimal Mass Transport} 
\author{Alistair Barton\thanks{This work is part of a Master's thesis prepared by A. Barton under the supervision of N. Ghoussoub.} \quad and \quad 
Nassif  Ghoussoub\thanks{Partially supported by a grant from the Natural Sciences and Engineering Research Council of Canada.}\\ \\
{\it\small Department of Mathematics,  University of British Columbia}\\
{\it\small Vancouver BC Canada V6T 1Z2}\\
}
\date{February 12, 2018}
\maketitle

\begin{abstract}  We investigate how mass transports that optimize the inner product cost -considered by Y. Brenier- propagate in time along a given Lagrangian. In the deterministic case, we consider 
 transports that maximize and minimize the following ``ballistic"  cost functional on phase space $M^*\times M$,  
\[
b_T(v, x):=\inf\{\langle v, \gamma (0)\rangle +\int_0^TL(t, \gamma (t), {\dot \gamma}(t))\, dt; \gamma \in C^1([0, T), M);   \gamma(T)=x\}, 
\]
where $M=\R^d$, $T>0$, and $L:M\times M \to \R$ is a suitable Lagrangian. We also consider the stochastic counterpart:
\begin{align*}
    \underline{B}_T^s(\mu,\nu):=\inf\left\{\mathbb{E}\left[\langle V,X_0\rangle +\int_0^T L(t, X,\beta(t,X))\,dt\right]; X\in \mathcal{A}, V\sim\mu,X_T\sim \nu\right\}
\end{align*}
where $\mathcal{A}$ is the set of stochastic processes satisfying $dX=\beta_X(t,X)\,dt+ dW_t,$ for some 
drift $\beta_X(t,X)$, and where $W_t$ is $\sigma(X_s:0\le s\le t)$-Brownian motion. While inf-convolution allows us to easily obtain Hopf-Lax formulas on Wasserstein space for cost minimizing transports, this is not the case for total cost maximizing transports, which actually are  sup-inf problems. However,  in the case where the Lagrangian $L$ is jointly convex on phase space,  Bolza-type dualities --well known in the deterministic case but novel in the stochastic case--transform sup-inf problems to sup-sup settings. Hopf-Lax formulas relate optimal ballistic transports to those associated with dynamic fixed-end transports studied by Bernard-Buffoni and Fathi-Figalli in the deterministic case, and by Mikami-Thieullen in the stochastic setting. 
We also write Eulerian formulations and point to links with the theory of mean field games.  

\end{abstract}


\section{Introduction and main results} Given a cost functional $c(y, x)$ on some product measure space $X_0\times X_1$, and  two probability measures $\mu$ on $X_0$ and $\nu$ on  $X_1$, we consider the problem of optimizing the total cost of {\it transport plans} and its corresponding dual principle as formulated by Kantorovich
\begin{equation*}
\inf\big\{\int_{X_0\times X_1} c(y, x)) \, d\pi; \pi\in \mK(\mu,\nu)\big\}=\sup\big\{\int_{X_1}\phi_1(x)\, d\nu(x)-\int_{X_0}\phi_0(y)\, d\mu(y);\,  \phi_1, \phi_0 \in \mK(c)\big\},
\end{equation*}
where $\mK(\mu,\nu)$ is the set of {\it transport plans} between $\mu$ and $\nu$, that is the set of probability measures $\pi$ on $X_0\times X_1$ whose marginal on $X_0$ (resp. on $X_1$) is $\mu$ (resp., $\nu$). On the other hand,  $\mK(c)$ is the set of functions $\phi_1\in L^1(X_1, \nu)$ and $\phi_0\in L^1(X_0, \mu)$ such that 
$
\phi_1(x)-\phi_0(y) \leq c(y,x)$ for  all $(y,x)\in X_0\times X_1.
$
The pairs of functions in $\mK(c)$ can be assumed to satisfy 
\begin{equation}\label{KD.condition}
\phi_1(x)=\inf_{y\in X_0} c(y, x)+\phi_0(y) \quad {\rm and} \quad\phi_0(y)=\sup_{x\in X_1} \phi_1(x)-c(y, x).
\end{equation}
They will be called {\it admissible Kantorovich potentials}, and for reasons that will become clear later, we shall say that $\phi_0$ (resp., $\phi_1$) is an initial (resp., final) Kantorovich potential.\\
The original Monge problem dealt with the cost $c(y, x)=|x-y|$ (\cite{M}, \cite{S}, \cite{E-G}, \cite{V1}, \cite{V2}) 
and was constrained to those probabilities in ${\mathcal K}(\mu, \nu)$ that are supported by graphs of measurable maps from $X$ to $Y$ pushing $\mu$ onto $\nu$.  Brenier \cite{B1} considered the important quadratic case  $c(x,y)=|x-y|^2$. This was followed by a large number of results addressing costs of the form $f(x-y)$, where $f$ is either a convex or a concave function \cite{G-M}. With a purpose of connecting mass transport with Mather theory, Bernard and Buffoni \cite{B-B} considered dynamic cost functions on a given compact manifold $M$, that deal with fixed end-points problems of the following type:
\begin{equation}\label{BB}
c_T(y,x):=\inf\{\int_0^TL(t, \gamma(t), {\dot \gamma}(t))\, dt; \gamma\in C^1([0, T), M);  \gamma(0)=y, \gamma(T)=x\},
\end{equation}
where $[0, T]$ is a fixed time interval, and $L: TM \to \R\cup\{+\infty\}$ is a given Lagrangian that is convex in the second variable of the tangent bundle $TM$. Fathi and Figalli \cite{F-F} eventually dealt with the case where $M$ is a non-compact Finsler manifold.  Note that standard cost functionals of the form  $f(|x-y|)$, where $f$ is convex, are particular cases of the dynamic formulation, since they correspond to Lagrangians of the form $L(t, x, p)=f(p)$. \\
We shall assume throughout that $M=M^*=\R^d$, while preserving --for pedagogical reasons-- the notational distinction between the state space and its dual. In this paper, we shall consider the {\it ``ballistic cost function,"} which is defined on phase space $M^*\times M$ by, 
 \begin{equation}\label{bal}
b_T(v, x):=\inf\{\langle v, \gamma (0)\rangle +\int_0^TL(t,  \gamma (t), {\dot \gamma}(t))\, dt; \gamma \in C^1([0, T), M);   \gamma(T)=x\}, 
\end{equation}
where $M$ is  a Banach space and $M^*$ is its dual. The associated transport problems will be 
\begin{equation}\label{bal2}
{\underline B}_T(\mu_0, \nu_T):=\inf \{\int_{M^*\times M} b_T(v, x)\, d\pi;\, \pi\in \mK(\mu_0,\nu_T)\}, 
\end{equation}
where $\mu_0$ (resp., $\nu_T$) is a given probability measure on $M^*$ (resp., $M$), 
and
\begin{equation}\label{bal1}
{\overline B}_T(\mu_0, \nu_T):=\sup \{\int_{M^*\times M} b_T(v, x)\, d\pi;\, \pi\in \mK(\mu_0,\nu_T)\}. 
\end{equation}
Note that when $T=0$, we have $b_0(x, v)=\langle v, x\rangle$, which is exactly the case considered by Brenier \cite{B1}, that is 
\begin{equation}
{\underline W}(\mu_0, \nu_0):=\inf \{\int_{M^*\times M} \langle v, x\rangle\, d\pi;\, \pi\in \mK(\mu_0,\nu_0)\}, 
\end{equation}
and 
\begin{equation}
{\overline W}(\mu_0, \nu_0):=\sup \{\int_{M^*\times M} \langle v, x\rangle\, d\pi;\, \pi\in \mK(\mu_0,\nu_0)\},  
\end{equation}
making (\ref{bal1}) a suitable dynamic version of the Wasserstein distance.\\
We shall also consider stochastic versions of the above problems, namely the cost of transport between two random variables $Y$ and $Z$ in $L^2(\Omega, M)$ defined as
\begin{equation}
c^s_T(Y, Z)=\inf\left\{\mathbb{E}\left[\int_0^T L(t,X,\beta(t,X))\,dt\right]; X\in \mathcal{A}, X_0=Y, X_T=Z\,  {\rm a.s}\right\},
\end{equation}
as well as the ballistic cost of using an input $V$ in $L^2(\Omega, M^*)$ to get to the random state $Z$ in $L^2(\Omega, M)$, namely
\begin{equation}
b^s_T(V, Z)=\inf\left\{\mathbb{E}\left[\langle V,X_0\rangle +\int_0^T L(t,X,\beta(t,X))\,dt\right]; X\in \mathcal{A}, X_T=Z\,  {\rm a.s}\right\},
\end{equation}
 where 
$\mathcal{A}$
is the set of stochastic processes verifying the stochastic differential equation 
$$dX=\beta_X(t,X)\,dt+ dW_t,$$ for some 
drift $\beta_X(t,X)$, where $W_t$ is $\sigma(X_s:0\le s\le t)$-Brownian motion. The corresponding mass transports are then 
\begin{align}\label{eq:stoc}
    C_T^s(\nu_0,\nu_T):=&\inf \left\{c^s_T(Y, Z); Y\sim\nu_0, Z\sim \nu_T\right\}\\
    =&\inf\left\{\expect{\int_0^T L(t,X,\beta_X(t,X))\,dt}; X\in \mathcal{A},\,  X_0\sim\nu_0,\, X_T\sim\nu_T\right\}, 
\end{align}
 which was considered by Mikami and Thieullen \cite{M-T}, while 
\begin{align}
    \underline{B}_T^s(\mu_0,\nu_T):&=\inf \left\{b^s_T(V, Z); V\sim\mu_0, Z\sim \nu_T\right\}\\
    &= \inf\left\{\mathbb{E}\left[\langle V,X_0\rangle +\int_0^T L(t,X,\beta(t,X))\,dt\right]; X\in \mathcal{A}, V\sim\mu_0,X_T\sim \nu_T\right\},\\ 
\overline{B}_T^s(\mu_0,\nu_T):&=\sup \left\{b^s_T(V, Z); V\sim\mu_0, Z\sim \nu_T\right\}\\
   &= \sup\limits_{V\sim\mu_0, Z\sim \nu_T}\inf\limits_{X\in \mathcal{A}, X_T=Z}\left\{\mathbb{E}\left[\langle V,X_0\rangle +\int_0^T L(t,X,\beta(t,X))\,dt\right]
  \right\},
\end{align}
that we shall consider in the sequel.\\
 In Section 2, we shall prove the following interpolation formulae on Wasserstein space associated to the deterministic minimization problem: 
\begin{equation}\label{HL.1}
{\underline B}_T(\mu_0,\nu_T)=\inf\{{\underline W}(\mu_0, \nu)+ C_T(\nu, \nu_T);\, \nu\in {\mathcal P}(M)\}.
\end{equation}
 The above  formula can be seen as extensions of those by Hopf-Lax on state space to Wasserstein space. Indeed, for any (initial) function $g$, the associated value function
 can be written as 
\begin{equation}\label{HL}
\phi_g(t,x)=\inf\{g(y)+c_t(y, x);\, y\in M\}. 
\end{equation}
In the case where the Lagrangian $L(t, x, p)=L_0(p)$ is only a function of $p$, and if $H_0$ is the associated Hamiltonian, then $c_t(y, x)=tL_0(\frac{1}{t}|x-y|)$
 and (\ref{HL}) is nothing but the Hopf-Lax formula used to generate solutions for corresponding Hamilton-Jacobi equations. When $g$ is the linear functional $g(x)=\langle v, x\rangle$, then $b_t( v, x)$ is itself a solution to the Hamilton-Jacobi equation, since 
\begin{equation}\label{basic}
b_t( v, x)=\inf\{\langle v, y\rangle+c_t(y, x);\, y\in M\}.
\end{equation}
In other words,  (\ref{HL.1})  
can now be seen as extensions of (\ref{basic}) to the space of probability measures, where the Wasserstein distance fill the role of the scalar product. \\
 In order to establish duality formulas, we consider the following forward Hamilton-Jacobi equations: \begin{eqnarray}\label{HJ+} 
\left\{ \begin{array}{lll}
\partial_t\phi+H(x, \nabla_x\phi)&=&0 \,\, {\rm on}\,\,  [0, T]\times M,\\
\hfill \phi(0, x)&=&f(x), 
\end{array}  \right.
 \end{eqnarray}
 and backward 
 Hamilton-Jacobi equations: 
 \begin{eqnarray}\label{HJ-} 
\left\{ \begin{array}{lll}
\partial_t\phi+H(x, \nabla_x\phi)&=&0 \,\, {\rm on}\,\,  [0, T]\times M,\\
\hfill \phi(T, x)&=&f(x), 
\end{array}  \right.
 \end{eqnarray}
where the Hamiltonian on $[0, T] \times M\times M^*$ is defined by
$
H(t, x, q)=\sup_{p\in M}\{\langle p, q\rangle -L(t,  x, p)\}.
$
Unless specified otherwise, we shall consider  ``variational solutions" for (\ref{HJ+}) and (\ref{HJ-}), which are formally given by the formulae 
\begin{equation}\label{value+}
\Phi^t_{_{f,+}}(x):=\Phi_{_{f,+}}(t,x)=\inf\Big\{f(\gamma (0))+\int_0^tL(s,\gamma (s), {\dot \gamma}(s))\, ds; \gamma \in C^1([0, T), M);   \gamma(t)=x\Big\}, 
\end{equation}
\begin{equation}\label{value-}
\Phi^t_{_{f,-}}(x):=\Phi_{_{f,-}}(t,x)=\sup\Big\{f(\gamma (T))-\int_t^TL(s,\gamma (s), {\dot \gamma}(s))\, ds; \gamma \in C^1([0, T), M);   \gamma(t)=x\Big\}. 
\end{equation}
Additional conditions on the Lagrangian are needed in order to verify if $\Phi_{_{f,+}}$ and $\Phi_{_{f,-}}$ are anywhere close to a classical solution. We shall then prove the following duality formulae:
 \begin{eqnarray}\label{ballistic.dual1}
{\underline B}_T(\mu_0,\nu_T)&=&\sup\left\{\int_M\Phi_{f_{_*}, +}(T,x)\, d\nu_T(x)+ \int_{M^*} {f}(v)\, d\mu_0(v); \,  \hbox{$f$ concave  in Lip(M$^*$)}
\right\}\\
&=&\sup\left\{\int_M g(x)\, d\nu_T(x)+\int_{M^*} (\Phi^0_{g,-})_{_*} (v)\,d\mu_0(v);\, \hbox{$g$ in ${\rm Lip(M)}$}\right\},
\end{eqnarray}
where ${h_*}$ is the concave Legendre transform of  $h$, i.e., 
$
{h_*}(v)=\inf\{\langle v, y\rangle - h(y); \, y\in M\}.$\\
 As to the question of attainment, we use a result by Fathi-Figalli \cite{F-F} to show that if $L$ is a Tonelli Lagrangian, and if $\mu_0$ is absolutely continuous with respect to Lebesgue measure, then there exists a probability measure $\pi_0$ on $M^*\times M$, and a concave function $k: M \to \R$ such that   
${\underline B}_T(\mu_0,\nu_T)=\int _{M^*} b_T \large(v, x) d\pi_0,$
and $\pi_0$ is supported on the possibly set-valued map $v\to \pi^*\phi^H_T(\nabla k_*(v), v)$, 
with $\pi^*:M\times M^*\to M$ being the canonical projection, 
and $(x, v) \to \phi^H_t(x, v)$ is the corresponding Hamiltonian flow.\\
In Section 3, we prove an analogous Hopf-Lax formulae on Wasserstein space associated to the stochastic minimization problem: 
\begin{equation}\label{SHL.one}
{\underline B}^s_T(\mu_0,\nu_T)=\inf\{{\underline W}(\mu_0, \nu)+ C^s_T(\nu, \nu_T);\, \nu\in {\mathcal P}(M)\}.
\end{equation}
As to the duality, there are two features that distinguish the deterministic case from the stochastic case. For one, there is no Monge-Kantorovich duality for the latter since it doesn't correspond to a cost minimizing transport problem. Moreover, stochastic processes are not reversible as deterministic paths and so we can only prove the following duality formula: 
\begin{equation}
   \lbstoc_T^s(\mu_0,\nu_T)=\sup\left\{\int_M g(x)\, d\nu_T(x)+\int_{M^*} (\Psi^0_{g,-})_{_*} (v)\,d\mu_0(v);\, \hbox{g in Lip(M)}\right\},
\end{equation}
 where  
this time 
$\Psi_{_{g, -}}$ is the solution to the backward Hamilton-Jacobi-Bellman equation (\ref{HJB}).
\begin{eqnarray}\label{HJB}
\left\{ \begin{array}{lll}
\partial_t\psi+\frac{1}{2}\Delta \psi+H(x, \nabla_x\psi)&=&0 \,\, {\rm on}\,\,  [0, T]\times M,\\
\hfill \psi(T, x)&=&g(x), 
\end{array}  \right.
 \end{eqnarray}
 whose formal variational solutions are given by the formula:
 \begin{equation}
    \Psi_{g,-}(t,x)=\sup_{X\in\mathcal{A}}\left\{\expcond{g(X(T))-\int_t^T L(s,X(s),\beta_X(s,X))\,ds}{X(t)=x}\right\}.
\end{equation}
In order to deal with the maximization problems ${\overline B}_T(\mu_0,\nu_T)$ and ${\overline B}^s_T(\mu_0,\nu_T)$, we need to use Bolza-type duality to convert the sup-inf problem to a concave maximization problem. For that, we shall assume that the Lagrangian $L$ is jointly convex in both variables. In Section 4, we then consider the dual Lagrangian ${\tilde L}$ defined on $M^*\times M^*$ by 
$$
\tilde L(t, v,q):=L^*(t, q, v)=\sup\{\langle v, y\rangle +\langle p, q\rangle -L(t, y, p);\, (y, p)\in M\times M\}, 
$$
and the corresponding fixed-end costs on $M^*\times M^*$, 
\begin{equation}\label{tilde}
{\tilde c}_T(u, v):=\inf\{\int_0^T{\tilde L}(t, \gamma (t), {\dot \gamma}(t))\, dt; \gamma\in C^1([0, T), M^*);  \gamma (0)=u, \gamma (T)=v\},
\end{equation}
and its associated transport
\begin{equation}\label{BBT*}
{\tilde C}_T(\mu_0, \mu_T):=\inf \{\int_{M^*\times M^*} {\tilde c}_T(x,y)\, d\pi;\, \pi\in \mK(\mu_0,\mu_T)\}. 
\end{equation}
We then recall the deterministic Bolza duality, and establish a new stochastic Bolza duality.\\
We use these results in Section 5, to establish the following results for ${\overline B}_T(\mu_0,\nu_T)$.  
\begin{equation}
{\overline B}_T(\mu_0,\nu_T)=\sup\{ {\overline W}(\nu_T, \mu)-{\tilde C}_T(\mu_0, \mu);\, \mu\in {\mathcal P}(M^*)\}.  
\end{equation}
and
\begin{equation}
{\overline B}_T(\mu_0,\nu_T)=\inf\left\{\int_Mg(x)\, d\nu_T(x)+ \int_{M^*} {\tilde\Phi}^0_{g^*,-}(v)\, d\mu_0(v); \,  \hbox{$g$ convex on $M$} 
\right\},
\end{equation}
where $g^*$ is the convex Legendre transform of  $g$, i.e., 
$
g^*(x)=\sup\{\langle v, x\rangle -  g(v);\, v\in M^*\}, 
$
and $\tilde\Phi_{k, -}$ is a solution 
 of the following  dual backward Hamilton-Jacobi equation:
\begin{eqnarray}\label{dHJ} 
\left\{ \begin{array}{lll}
\partial_t\phi-H(t,\nabla_v\phi, v)&=&0 \,\, {\rm on}\,\, [0, T]\times M^*,\\
\hfill \phi(T, v)&=&k(v),  
\end{array}\right.
\end{eqnarray} 
whose variational solution is given by
\begin{equation}\label{value.2}
\tilde \Phi_{k,-}(t,v)=\sup\Big\{k(\gamma (T))-\int_0^t{\tilde L}(s,  \gamma (s), {\dot \gamma}(s))\, ds; \gamma \in C^1([0, T), M^*);   \gamma(0)=v\Big\}. 
\end{equation}
In Section 6, we deal with the stochastic counterpart ${\overline B}^s_T(\mu_0,\nu_T)$ and prove the following 
\begin{equation}\label{SHL.three}
   \overline{B}_T^s(\mu_0,\nu_T):=\sup\left\{\mathbb{E}\left[\langle X,V(T)\rangle -\int_0^T {\tilde L}(t,V,\beta(t,V))\,dt\right]; V\in {\mathcal A}, V_0\sim\mu_0,X \sim \nu_T\right\}
\end{equation}
and therefore
\begin{equation}\label{SHL.two}
{\overline B}^s_T(\mu_0,\nu_T)=\sup\{ {\overline W}(\nu_T, \mu)-{\tilde C}^s_T(\mu_0, \mu);\, \mu\in {\mathcal P}(M^*)\}, 
\end{equation}
as well as the following duality formula:
\begin{equation}
    \overline{B}_T^s(\mu_0,\mu_T)=\inf\left\{\int_{M^*} g(x)\,d\nu_T+\int_M{\tilde\Psi}^0_{g^*,-}(v)\,d\mu_0;\, \, \hbox{$g$ convex in $\dbound (M^*)$} \right\},
\end{equation}
where $\tilde\Psi_k$ solves the Hamilton-Jacobi-Bellman equation
\begin{eqnarray}\label{HJB2} 
\left\{ \begin{array}{lll}
\partial_t\psi+\frac{1}{2}\Delta \psi-H(\nabla_v\psi, v)&=&0 \,\, {\rm on}\,\,  [0, T]\times M^*,\\
\hfill \psi(T, v)&=&k(v), 
\end{array}  \right.
 \end{eqnarray}
 whose formal variational solutions are given by the formula:
 \begin{equation}
 \label{eq:max.dynprog}
    \tilde\Psi_{k, -}(t,v)=\sup_{X\in\mathcal{A}}\left\{\expcond{k(X(T))-\int_t^T {\tilde L}(s,X(s),\beta_X(s,X))\,ds}{X(t)=v}\right\}.
\end{equation}
Finally, a few words about our notation: We shall denote by $\partial g$ the subdifferential of a convex function $g$, and by $\tilde\partial h:=-\partial (-h)$ the superdifferential of a concave function $h$.\\
 The set of probability measures on a Banach space $X$ will be denoted $\mathcal{P}(X)$, while the subset of those with finite first moment will be denoted 
$$\mathcal{P}_1(X):=\{\nu\in\mathcal{P}(X); \int_X\abs{x}\,d\nu(x)<\infty\}.$$
 $\mathcal{P}_1(X)$ is clearly a subset  of the Banach space of all finite measures with finite first moment, denoted similarly $\mathcal{M}_1(X):=\{\nu;\int_X 1+\abs{x}\,d\nu(x)<\infty\}$, which is dual to the Banach space  $\lip{X}$ of all bounded uniformly Lipschitz functions on $X$. 
For the stochastic part, we shall also need to work with the space $\dbound(X):=\lip{X}\cap C^\infty(X)$.

Several of the above results appeared in the posted but non-published manuscripts \cite{G}, which dealt with the deterministic case and \cite{B-G}, which addressed the stochastic case. We eventually elected to combine them in a single publication so as to illustrate the obvious similarities, but also the subtle differences between the two cases.

\section{Minimizing the ballistic cost: Deterministic case}
In this section we deal with the standard transportation problem associated to the cost $b_T(v,x)$. 
We shall assume that the Lagrangian $L$ satisfies the following:
\begin{trivlist}
 
\item $(A0)$ The Lagrangian $(t,x,v)\mapsto L(t,x,v)$ is bounded below, and for all $(t,x)\in [0,T]\times M$, $v\mapsto L(t,x,v)$ is convex and $\delta$-coercive in the sense that there is a $\delta>1$ such that 
\begin{equation}
    \lim_{\abs{v}\rightarrow \infty}\frac{L(t,x,v)}{\abs{v}^\delta}=+\infty.
\end{equation}
 
\end{trivlist}

\begin{thm} \label{duality} Assume that $L$ satisfies $(A0)$  and let $\mu_0$ (resp. $\nu_T$) be a probability measure on $M^*$  (resp., $M$) with finite first moment. Then, the following interpolation formula holds:
\begin{equation}\label{HL.one}
{\underline B}_T(\mu_0,\nu_T)=\inf\{{\underline W}(\mu_0, \nu)+ C_T(\nu, \nu_T);\, \nu\in {\mathcal P}_1(M)\}.
\end{equation}
The infimum is attained at some probability measure $\nu_0$ on $M$, and the initial Kantorovich potential for $C_T(\nu_0, \nu_T)$ is concave. 
 \end{thm}
{\bf Proof:} To prove the formula it suffices to note that 
\begin{eqnarray*}
   && \inf\left\{\underline{W}(\mu_0,\nu)+C_T(\nu,\nu_T);\, \nu\in {\mathcal P}_1(M)\right\}\\
    &&\qquad \qquad =\inf_{\nu\in {\mathcal P}_1(M)}\left\{\int_{M^*\times M}\langle v,x\rangle\,d\pi_W(v,x) +\int_{M\times M} c_T(x,y)\,d\pi_C(x,y);\pi_W\in\mathcal{K}(\mu_0,\nu),\pi_C\in\mathcal{K}(\nu,\nu_T)\right\}\\
   && \qquad \qquad =\inf\limits_{\pi\in {\mathcal P}_1(M^*\times M\times M)}\left\{\int_{M^*\times M\times M}\langle v,x\rangle + c_T(x,y)\,d\pi(v,x,y);\pi_{1}=\mu_0,\pi_{3}=\nu_T\}\right\}\\
 && \qquad \qquad \ge  \underline{B}(\mu_0,\nu_T).
\end{eqnarray*}
For the reverse inequality, use your favourite selection theorem to find a measurable function $y_\epsilon:M^\ast\times M\rightarrow M$ that satisfies $\langle v,y_\epsilon(v,x)\rangle +c(y_\epsilon(v,x),x)-\epsilon<b_T(v,x)$. Fixing $\pi\in\mathcal{K}(\mu_0,\nu_T)$ and letting $\pi_\epsilon:=(\text{Id}\times\text{Id}\times y_\epsilon)_\# \pi \in\mathcal{P}(M^\ast\times M\times M)$
\begin{equation*}\begin{split}
    \underline{B}(\mu_0,\nu_T)=&\int_{M^*\times M} b_T(v,x)\,d\pi(v,x)\ge \int_{M^*\times M\times M} \langle v,y\rangle +c_T(y,x)\,d\pi_\epsilon(v,x,y)-\epsilon.
\end{split}
\end{equation*}
To show that the minimizer is achieved, we need to prove that $C_T$ satisfies a coercivity condition on the space $\mathcal{P}_1(M)$ of probabilities on $M$ with finite first moments. 
For that, we show that for any fixed $\nu_T\in\mathcal{P}_1(M)$ and any positive constant $N>0$, the set of measures $\nu\in \mathcal{P}_1(M)$ satisfying 
\begin{equation}\label{eq:det.tight}
   C_T(\nu,\nu_T)\le N\int_M \abs{x}\,d\nu(x),
\end{equation}
is tight. Indeed, 
from $(A0)$,  there exists a constant $K$ such that $c_T(x,y)>N\abs{\frac{x-y}{T}}^\delta-K$. We concern ourselves with the cylinder set $B:=\ball{R}{0}^c\times M$. Let $\nu\in\mathcal{T}_{\epsilon,R}:=\{\nu;\nu(\ball{R}{0}^c)>\epsilon\}$. We shall assume, without loss of generality, that $L$ and hence $c_T$ is non-negative, hence 
for any optimal transport plan $\pi\in \mathcal{K}(\nu,\nu_T)$
\begin{equation}\begin{split}\label{det.Coercivity.eq}
    C_T(\nu,\nu_T)\ge&  \frac{N}{T^\delta}\int_{B}\abs{\abs{x}-\abs{y}}^\delta\,d\pi(x,y)-K\epsilon.
\end{split}
\end{equation}
We define $\hat{\pi}:=\pi\rvert_B$ to be the restriction of $\pi$ to the set $B$, and transfer the problem 
 to $\mathbb{R}_+$ by using the push-forward $\bar{\pi}:=(\abs{\cdot}\times\abs{\cdot})_\#\hat{\pi}$ to obtain,
\begin{equation}\begin{split}
    \int_{B}\abs{\abs{x}-\abs{y}}^\delta\,d\hat{\pi}(x,y)=\int_{\mathbb{R}_+\times \mathbb{R}_+} \abs{x-y}^\delta\,d\bar{\pi}(x,y).
\end{split}
\end{equation}
 We can obtain a lower estimate for this by minimizing over transportation measures sharing $\bar{\pi}$'s marginals (i.e., $\gamma\in \mathcal{K}(\bar{\pi}_1,\bar{\pi}_2)$). This is a well known optimal transport problem, whose optimal plan given by the monotone Hoeffding-Frechet mapping $x\mapsto G_{\nu}(G_{\nu_T}^{-1}(x))$, where $G_\nu(t):=\inf\{z\in\mathbb{R}:t\ge\nu(\{x\le z\})\}$ is the quantile function associated with the measure $\nu$ \cite{beiglbock}. Thus the optimal plan maps each quantile in one measure to the corresponding quantile in the other.
Substituting this into the integral and applying Jensen's inequality:
\begin{equation}\begin{split}\label{eq:det.coerce}
    \int_{\mathbb{R}_+\times \mathbb{R}_+} \abs{x-y}^\delta\,d\bar{\pi}(x,y)\ge & \int_{\mathbb{R}_+\times \mathbb{R}_+}\abs{x-y}^\delta\,d\bracket{\bracket{G_{\abs{Y(0)}}\times G_{\abs{Y(T)}} }_\#\lambda_{[0,1]}}(x,y)\\
    \ge&  \bracket{\int_{\ball{R}{0}^c}\abs{x}\,d\nu(x)-b(\nu_T)}^\delta,
\end{split}
\end{equation}
where $b(\nu_T):=\int\abs{x}\,d\nu_T$ and $R>b/\epsilon$. We thus want to find $R$ such that
\begin{equation*}
    N \bracket{\bracket{\int_{\ball{R}{0}^c}\abs{x}\,d\nu(x)-b(\nu_T)}^\delta-K\epsilon}>N\bracket{\int_{\ball{R}{0}^c}\abs{x}\,d\nu(x)+R(1-\epsilon)}\ge N\int\abs{x}\,d\nu(x).
\end{equation*}
Letting $I_\nu(R):=\int_{\ball{R}{0}^c}\abs{x}\,d\nu(x)$ ($\ge R\epsilon$ for $\nu\in\mathcal{T}_{\epsilon,R}$), we find the condition
\begin{equation*}
    \frac{\alpha}{T^{\delta-1}}  \bracket{\bracket{I_\nu(R)^{1-\frac{1}{\delta}}-b(\nu_T)I_\nu(R)^{-\frac{1}{\delta}}}^\delta-\frac{\bracket{UT}}{I_\nu(R)}}>N\bracket{1+\frac{R(1-\epsilon)}{I_\nu(R)}},
\end{equation*}
which by using the mentioned bound on $I_\nu(R)$ is satisfied if $R$ is large enough so that 
\begin{equation*}
    \frac{\alpha}{T^{\delta-1}} \bracket{\bracket{(R\epsilon)^{1-\frac{1}{\delta}}-b(\nu_T)(R\epsilon)^{-\frac{1}{\delta}}}-\frac{\bracket{UT}}{R\epsilon}}>\frac{N}{\epsilon}.
\end{equation*}
To show the minimizer is achieved, fix any $\nu_0$ in $\mathcal{P}_1(M)$, and note that by coercivity the set of probability measures $\nu$ such that
\begin{equation*}
    \underline{W}(\mu_0,\nu)+C_T(\nu,\nu_T)\le b(\mu_0)\int_M\abs{x}\,d\nu+C_T(\nu,\nu_T)\leq b(\mu_0)\int_M\abs{x}\,d\nu+C_T(\nu_0,\nu_T)
\end{equation*}
is tight.
\begin{rmk}\rm
Note that (\ref{eq:det.coerce})
indicates that when $\nu_1\in\mathcal{P}_1(M)$ and $\nu_0\in\mathcal{P}(M)\setminus\mathcal{P}_1(M)$, then $C(\nu_0,\nu_1)=C(\nu_1,\nu_0)=\infty$. 
\end{rmk}

\begin{thm} \label{duality} Assume that $L$ satisfies $(A0)$ and let $\mu_0$ (resp. $\nu_T$) be a probability measure on $M^*$  (resp., $M$) with finite first moment.
\begin{enumerate}

\item  If $\mu_0$ has compact support, then we have the following duality formula
\begin{equation}
{\underline B}_T(\mu_0,\nu_T)=\sup\left\{\int_M g(x)\, d\nu_T(x)+\int_{M^*} (\Phi^0_{g,-})_*(v)\,d\mu_0(v);\, \hbox{g in Lip(M$$)}\right\}.
\end{equation}
 \item  If $\nu_T$ has compact support, then 
\begin{equation}
{\underline B}_T(\mu_0,\nu_T)=\sup\left\{\int_M\Phi^T_{f_*, +}(x)\, d\nu_T(x)+ \int_{M^*} {f}(v)\, d\mu_0(v); \,  \hbox{$f$ concave  in Lip(M$^*$)}
\right\}.
\end{equation}

\end{enumerate}
\end{thm}
\noindent {\bf Proof:} We shall need the following identifications of the Legendre transforms in the Banach space  $\mathcal{M}_1(\R^n)$  of measures $\nu$ on $\R^n$ such that $\int_{\R^n} \bracket{1+\abs{x}}\,d\nu<\infty$ in duality with the space of Lipschitz functions $\lip{\R^n}$. 

 \begin{lem}\label{lem:Wlegendre}
a) For $\mu_0\in\mathcal{P}(\R^n)$ with compact support, define $\underline{W}_{\mu_0}:\mathcal{M}_1(\R^n)\rightarrow \mathbb{R}\cup\{\infty\}$  
to be
\begin{equation*}
    \underline{W}_{\mu_0}(\nu):=\begin{cases} \underline{W}(\mu_0,\nu)&\nu\in\mathcal{P}_1(M)\\
   + \infty&\text{otherwise}.
    \end{cases}
\end{equation*}Then,  the convex Legendre transform of $\underline{W}_{\mu_0}$ is given for $f\in \lip{\R^n}$ by 
 $   \underline{W}_{\mu_0}^\ast(f)=-\int_M f_*\,d\mu_0.$
 
 b) For $\nu_0\in\mathcal{P}(\R^n)$, define the function $C_{\nu_0}:\mathcal{M}_1(\R^n)\rightarrow\mathbb{R}\cup\{\infty\}$ to be
\begin{equation}
    C_{\nu_0}(\nu):=\begin{cases}C_T(\nu_0,\nu)&\nu\in\mathcal{P}_1(\R^n)\\
    +\infty&\text{otherwise.}\end{cases}
\end{equation}
Then,  the convex Legendre transform of $C_{\nu_0}$ is given for $f\in \lip{\R^n}$ by  $   C_{\nu_0}^\ast(f)=\int_M \phi_{_{f,-}}(0, x)\,d\nu_0(x),$
where $\phi_{_{f,-}}$ is the solution to the backward Hamilton-Jacobi equation (\ref{HJ-}) with final condition $\phi_{_{f,-}}(T,x)=f(x)$. 

\end{lem}
\noindent{\bf Proof:} Both statements follow from Kantorovich duality. Indeed, both functions are convex and weak$^*$-lower semi-continuous on $\mathcal{M}_1(\R^n)$. Since $\mu_0$ has compact support, Brenier's duality yields 
\begin{equation*}
    \underline{W}_{\mu_0}(\nu)=\sup_{g\in\lip{M}}\left\{\int_M g\,d\nu +\int_{M^*} g_*\,d\mu_{0}\right\}.
\end{equation*}
We then have
\begin{equation}
   \underline{W}_{\mu_0}^\ast(f)=\sup_{\nu\in\mathcal{M}_1(M)}\inf_{g\in\lip{M}}\{\int_M f\,d\nu -\int_M g\,d\nu-\int_{M^*} g_*\,d\mu_{0}\}.
\end{equation}
Note that the functional $g\mapsto -\int_{M^*} g_{_{*}} \,d\mu=\int_{M^*} (-g)^\ast \,d\hat{\mu}(v)$ (where $d\hat{\mu}(v):=d\mu(-v)$) is convex and lower semicontinuous,
and we may therefore apply the Von Neuman minimax theorem as the expression is linear in $\nu$ and convex in $g$. We obtain
\begin{equation}
    \underline{W}_{\mu_0}^\ast(f)=\inf_{g\in\lip{M}}\sup_{\nu\in\mathcal{M}_1(M)}\{\int_M f\,d\nu -\int_M g\,d\nu-\int_{M^*} g_*\,d\mu_{0}\}.
\end{equation}
The infimum must occur at $g=f$ since otherwise the sup in $\nu$ is $+\infty$, resulting in  statement a).

The same proof applies to $ C_{\nu_0}$, since in view of  the duality formula of Bernard and Buffoni  \cite{B-B}[Proposition 21]:
\begin{equation}
    C_{\nu_0}(\nu)=\sup_{g\in\lip{M}}\{\int_M g\,d\nu-\int_M \Phi^0_{_{g,-}}\,d\nu_0\}.
\end{equation}
Note that this holds for all $\nu\in \mathcal{M}_1(M)$, since if $g$ solves HJ, then so does $g+c$ for arbitrarily large $c$. 
We may again apply the minimax theorem as the expression is linear in $\nu$ and convex in $g$. To complete the proof of the theorem, we first note that Kantorovich duality yields that 
     $\nu\mapsto\lbstoc(\mu_0,\nu)$ is weak$^*$-lower semi-continuous on $\mathcal{P}_1(M)$ for all $\mu_0\in\mathcal{P}_1(M^\ast)$ and that $(\mu_0,\nu_T)\mapsto \lbstoc(\mu_0,\nu_T)$ is jointly convex.
Let now $ \underline{B}_{\mu_0}(\nu):=\underline{B}_T(\mu_0,\nu)$ if $\nu \in \mathcal{P}_1(M)$ and $+\infty$ otherwise. It follows that
\begin{equation}
    \underline{B}_{\mu_0}(\nu)=\ \underline{B}_{\mu_0}^{\ast\ast}(\nu):=\sup\{\int_{\R^n} f\,d\nu -\underline{B}_{\mu_0}^{\ast}(f);\, f\in \lip{\R^n}\}.
\end{equation}
 Now use the Hopf-Lax formula established above to write
\begin{equation}
\begin{split}
    \underline{B}_{\mu_0}^{\ast}(f):=&\sup\{\int_M f\,d\nu -\underline{B}_{\mu_0}(\nu); \nu \in \mathcal{P}_1(M)\}\\
    =&\sup\{\int_M f\,d\nu -\underline{W}(\mu_0,\nu')-C_T(\nu',\nu); \nu,\nu' \in \mathcal{P}_1(M) \}\\
     =&\sup\{\int_M \Phi^T_{_{f,-}}\,d\nu' -\underline{W}(\mu_0,\nu'); \nu' \in \mathcal{P}_1(M) \}\\
    =&-\int_M (\Phi^T_{_{f,-}})_*\,d\mu_0.
\end{split}
\end{equation}
This completes the proof of the first duality formula. \\
The second follows in the same way by simply varying the initial measure as opposed to the final measure in $\underline{B}_T(\mu, \nu)$. The concavity of $f$ follows from the Kantorovich dual condition (\ref{KD.condition}) and the linearity of $b_T$ in $v$.
\par

We now consider the problem of attainment for $\underline{B}_T(\mu, \nu)$. For that, 
we shall consider Tonelli Lagrangians studied in the compact case by Bernard-Buffoni \cite{B-B}, and by Fathi-Figalli \cite{F-F} in the case of a Finsler manifold. 

\begin{defn} \rm We shall say that $L$ is a {\it Tonelli Lagrangian on $M\times M$}, if it is $C^2$ and  satisfies (A0) with the additional requirement that the  function $v\to L(x, v)$ is strictly convex on $M$. 
\end{defn}
We also recall the following \cite[Definition 5.5.1, page 129]{A-G-S}:
\begin{defn} 
\label{approxdiff}
{\rm
Say that $f :M \rightarrow \R$ has an \textit{approximate differential} at $x \in M$
if there exists a function $h:M \rightarrow \R$
differentiable at $x$ such that the set $\{f = h\}$ has density $1$ at $x$ with respect to the Lebesgue measure.
In this case, the approximate value of $f$ at $x$ is defined as $\tilde f(x)=h(x)$, and the approximate differential of $f$ at $x$ is defined as $\tilde d_xf=d_xh$.
  It is not difficult to show that this definition makes sense. In fact, both $h(x)$, and $d_xh$ do not depend on the choice of $h$, provided $x$ is a density point of the set $\{f = h\}$.
}
\end{defn}
If $L$ is a Tonelli Lagrangian,
 the Hamiltonian $ H: M\times M^*\to\R$ is then C$^1$, and  the Hamiltonian vector field $X_H$ on $M\times M^*$ is then 
$
X_H(x,v)=(\frac{\partial H}{\partial v}(x,v),-\frac{\partial H}{\partial x}(x,v)), 
$
and the associated system of ODEs is given by
\begin{equation}\label{Ham}
\left\{
\begin{array}{l}
\dot x=\dfrac{\partial H}{\partial v}(x,v)\\
\dot v=-\dfrac{\partial H}{\partial x}(x,v).
\end{array}
\right.
\end{equation}
The connection between minimizers $\gamma:[a,b]\to M$ of $I_L$ and solutions of  $(\ref{Ham})$ is as follows. If we write 
$
x(t)=\gamma(t)\quad \text{and} \quad
v(t)=\frac{\partial L}{\partial p}(\gamma(t),\dot\gamma(t)), 
$
then 
$x(t)=\gamma(t)$ and $v(t)$ are C$^1$ with $\dot x(t)=\dot\gamma(t)$, and the Euler-Lagrange equation yields 
$
\dot v(t)=\frac{\partial L}{\partial x}(\gamma(t),\dot\gamma(t)),
$
from which follows 
that $t\mapsto (x(t),v(t))$ satisfies (\ref{Ham}). Note also that since $L$ is a Tonelli Lagrangian, the Hamiltonian $H$ is actually C$^2$, and  the vector field $X_H$ is C$^1$. It therefore defines a (partial) C$^1$ flow $\phi^H_t$. 

There is also a (partial) {\rm C}$^1$  flow $\phi^L_t$ on $M\times M^*$ such that every speed curve of an $L$-minimizer is a part of an orbit of $\phi^L_t$. This flow is called the Euler-Lagrange flow, is defined by
$\phi^L_t=\Leg^{-1}\circ \phi^H_t \circ \Leg,$
where $
\Leg: M\times M \rightarrow M\times M^*,
$
is the global Legendre transform
$
(x,p) \mapsto (x,\frac{\partial L}{\partial p}(x,p)).
$
Note that $\Leg$ is a homeomorphism on its image whenever $L$ is a Tonelli Lagrangian. 

\begin{thm} In addition to $(A0)$, assume that $L$ is a Tonelli Lagrangian and that $\mu_0$ is absolutely continuous with respect to Lebesgue measure. Then, 
there exists a concave function $k: M \to \R$ such that   
\begin{equation}
{\underline B}_T(\mu_0,\nu_T)=\int _{M^*} b_T \large(v, S_T \circ \nabla k_*(v)\large) d\mu_0(v), 
\end{equation}
where 
$
S_T(y)=\pi^*\phi^H_T(y,\nabla k(y)),
$
$\pi^*:M\times M^*\to M$ being the canonical projection,  and $\phi^H_t$ the Hamiltonian flow  associated to $L$. In other words, an optimal map for ${\underline B}_T(\mu_0,\nu_T)$ is given by $v\to \pi^*\phi^H_T(\nabla k_*(v), v)$. 
\end{thm}
\noindent{\bf Proof:} Start again by the interpolation inequality, 
$
{\underline B}_T(\mu_0,\nu_T)=C_T(\nu_0, \nu_T)+{\underline W}(\mu_0, \nu_0)
$
for some probability measure $\nu_0$. By the above and Kantorovich duality, there exists a concave function $k: M \to \R$ and another function $h:M\to \R$ such that $(\nabla  k_*)_\#\mu_0=\nu_0$,
\[
{\underline W}(\mu_0, \nu_0)=\int_M \langle \nabla k_*(v), v\rangle d\mu_0(v), 
\]
and
\[
C_T(\nu_0, \nu_T)=\int_Mh(x)\, d\nu_T(x)-\int_M k(y)\, d\nu_0(y).
\]
Now use a result of Fathi-Figalli \cite{F-F} to write 
$C_T(\nu_0, \nu_T)=\int _{M} c_T (y, S_Ty) d\nu_0(y),$ 
where $S_T(y)=\pi^*\phi^H_T(y,\tilde d_y k)$. 
Note that 
\begin{equation}
\hbox{${\underline B}_T(\mu_0,\nu_T)\leq \int _{M^*} b_T (v, S_T\circ \nabla k_*(v)) d\mu_0(v),$}
\end{equation}
since $(\nabla k_*)_\#\mu_0=\nu_0$ and $(S_T)_\#\nu_0=\nu_T$, and  therefore $(I \times S_T\circ \nabla k_*)_\#\mu_0$ belongs to ${\mathcal K}(\mu_0, \nu_T)$.\\
On the other hand, since $ b_T(v, x)\leq  c_T(\nabla k_*(v), x )+\langle \nabla k_*(v), v\rangle$ for every $v\in M^*$, we have
\begin{eqnarray*}
{\underline B}_T(\mu_0,\nu_T)&\leq& \int _{M^*} b_T (v, S_T\circ \nabla k_* (v)) d\mu_0(v)\\
&\leq& \int_{M^*}\{c_T(\nabla k_*(v), S_T\circ \nabla k_* (v)) +\langle \nabla k_*(v), v\rangle \}\, d \mu_0(v)\\
&=&\int _{M} c_T (y, S_Ty) d\nu_0(y) + \int_{M^*}\langle \nabla k_* (v), v\rangle \, d \mu_0(v)   \\
&=& C_T(\nu_0, \nu_T)+{\underline W}(\mu_0, \nu_0)\\
&=&{\underline B}_T(\mu_0,\nu_T).
\end{eqnarray*}
It follows that 
\begin{equation*}
{\underline B}_T(\mu_0,\nu_T)=\int _{M^*} b_T (v, S_T\circ \nabla k_* (v)) d\mu_0(v)=\int _{M^*} b_T (v, \pi^*\phi^H_T(\nabla k_* (v),\tilde d_{\nabla k_* (v)}k d\mu_0(v). 
\end{equation*}
Since $k$ is concave, we have that $\tilde d_xk=\nabla k(x)$, hence $\tilde d_{\nabla k_*(v)}k=\nabla k\circ \nabla k_*(v)=v$, which yields our claim that 
$
 {\underline B}_T(\mu_0,\nu_T)=\int _{M^*} b_T \large(v, \pi^*\phi^H_T(\nabla k_* (v), v)\large) d\mu_0(v).
 $

\section{Minimizing the ballistic cost: Stochastic case}

We now turn to the stochastic version of the minimizing cost. The methods of proof are generally similar to those for the deterministic cost, however there are two complications: The first is that stochastic mass transport does not fit in the framework of cost minimizing transports, hence the Kantorovich duality is not readily available. The second is that stochastic processes are not reversible and therefore there is only one direction to the transport, hence only one duality formula. In order to deal with the first complication, we rely on the results of Mikami-Thieullen  \cite{M-T} and therefore use the same assumptions that they imposed on the Lagrangian, namely

\begin{trivlist}

\item [(A1)] $L(t, x,v)$ is continuous, convex in $v$, and uniformly bounded below by a convex function $\underline{L}(v)$ that is $2$-coercive in the sense that
   $ \lim_{\abs{v}\rightarrow \infty}\frac{\underline{L}(v)}{\abs{v}^2}>0.$

\item [(A2)] $(t,x)\mapsto\log(1+L(x,u))$ is uniformly continuous in that
    \begin{equation*}
        \Delta L(\epsilon_1,\epsilon_2):=\sup_{u\in M^\ast}\left\{\frac{1+L(x,u)}{1+L(y,u)}-1;\abs{t-s}<\epsilon_1,\abs{x-y}<\epsilon_2\right\}\overset{\epsilon_1,\epsilon_2\rightarrow0}{\longrightarrow} 0.
    \end{equation*}
    \item [(A3)] The following boundedness conditions: 
    \begin{trivlist}
     \item (i) $\sup_{t,x} L(t,x,0)<\infty$.
     \item (ii) $\abs{\nabla_x L(t,x,v)}/(1+L(t,x,v))$ is bounded.
     \item (iii) $\sup\left\{\abs{\nabla_v L(t,x,u)}: \abs{u}\le R\right\}<\infty$ for all $R$.
    \end{trivlist}

\end{trivlist}
We will use the notation $X=(X_0,\beta_X,\sigma_X)$ to refer to an It\^o process $X(t)$ of the form:
\begin{equation}
    X(t)=X_0+\int_0^t \beta_X(s)\,ds +\int_0^t\sigma_X(s)\,dW_s.
\end{equation}We will use the notation $\mathcal{A}_{\nu_0}^{\nu_T}$ to refer to the set of stochastic processes $X=(X_0,\beta_X,\text{Id})$ with $X(0)\sim\nu_0$ and $X(T)\sim\nu_T$. Notably, (A1) implies that $\expect{L(t,X(t),\beta_X)}=\infty$ if $\beta_X(t)\not\in L^2(\mathbb{P})$.\\

Our main result for this section is the stochastic counterpart to Theorem \ref{duality}:  
\begin{thm}
\label{thm:stoch.dual} If $L$ satisfies the assumptions (A1), (A2), and (A3), then  
\begin{enumerate}
\item For any given probabilities $\mu_0\in\mathcal{P}(M^\ast)$ and $\nu_T\in\mathcal{P}(M)$, we have: 
\begin{equation}\label{eq:interpol}
    \lbstoc_T^s(\mu_0,\nu_T)=\inf \{\underline{W}(\mu_0,\nu)+C_T^s(\nu,\nu_T);\, \nu \in\mathcal{P}_1(M)\}.
    \end{equation}
Furthermore, this infimum is attained whenever $\mu_0\in\mathcal{P}_1(M^\ast)$ and $\nu_T\in\mathcal{P}_1(M)$.

\item If $\nu_T\in\mathcal{P}_1(M)$ and $\mu_0\in\mathcal{P}_1(M^\ast)$ are such that $\lbstoc(\mu_0,\nu_T)<\infty$, and if $\mu_0\in\mathcal{P}_1(M^\ast)$ has compact support, then  
\begin{equation}\label{eq:mindual}
    \lbstoc_T^s(\mu_0,\nu_T)=\sup\left\{\int_M f(x)\, d\nu_T(x)+\int_{M^*} (\Psi^0_{f,-})_*(v)\,d\mu_0(v);\, \hbox{$f\in \lip{M}$}\right\},
\end{equation}
where 
$\Psi_{_{f,-}}$ is the solution to the Hamilton-Jacobi-Bellman equation
\begin{align}\label{eq:HJB}\tag{HJB}
    \pderiv{\psi}{t}+\frac{1}{2}\Delta\psi(t,x)+H(t,x,\nabla\psi)=0, \qquad \psi(T,x)=&f(x).
\end{align}
\end{enumerate}
\end{thm}
{\bf Proof:}  
1) \,  First, expand $\underline{W}(\mu_0,\nu)$ and $C_T(\nu,\nu_T)$ in the interpolation formula to obtain:
\begin{eqnarray*}
   && \inf \{\underline{W}(\mu_0,\nu)+C_T^s(\nu,\nu_T);\, \nu \in\mathcal{P}_1(M)\}\\
&&\qquad \qquad  =\inf\left\{ \expect{\langle V,X\rangle +\int_0^T L(t,X(t),\beta_X(t))\,dt};V\sim \mu_0, X\sim\nu, X(t)\in\mathcal{A}_{\nu}^{\nu_T};\nu\in\mathcal{P}_1(M)
    \right\}\\
    &&\qquad \qquad  \leq \underline{B}(\mu_0,\nu_T).
    \end{eqnarray*}
To obtain the reverse inequality, let $\nu_n$ be a sequence of measures approximating the infimum in (\ref{eq:interpol}). Then for each $\nu_n$, there exists a stochastic process $Z_n\in\mathcal{A}_{\nu_n}^{\nu_T}$ such that
\begin{equation}
\expect{\int_0^T L(t,Z_n(t),\beta_{Z_n}(Z_n,t))\,dt}< C^s_T(\nu_n,\nu_T)+\tfrac{1}{n}.
\end{equation}
Similarly, let $d\gamma_x^n(v)\otimes d\nu_n(x)=d\gamma_n(v,x)$ be the disintegration of a measure $\gamma_n$ such that
\begin{equation*}
    \int \langle v,x\rangle \,d\gamma_n(v,x)<\underline{W}(\mu_0,\nu_n)+\tfrac{1}{n},
\end{equation*}
and define $U_n:M\times\Omega\rightarrow M^\ast$ to be a random variable such that $U_n[x]\sim \gamma_{x}^n$ for $\nu_n$-a.a. $x$. Thus $(U_n[Z_n(0)],Z_n(0))\sim \gamma_n$ and we have constructed a random variable that approximates the interpolation, as
\begin{equation}
\expect{\langle U[Z_n(0)],Z_n(0)\rangle +\int_0^T L(t,Z_n(t),\beta_{Z_n}(t))\,dt}\le \inf \{\underline{W}(\mu_0,\nu)+C_T^s(\nu,\nu_T);\, \nu \in\mathcal{P}_1(M)\}+\frac{3}{n}.
\end{equation}
To show that the infimum in $\nu$ is attained in the set $\mathcal{P}_1(M)$, we need again to prove the following coercivity property.\\
{\bf Claim:} For any fixed $\nu_T\in\mathcal{P}_1(M)$, $N\in\mathbb{R}$, the set of measures $\nu\in \mathcal{P}_1(M)$ satisfying 
  $  C(\nu,\nu_T)\le N\int \abs{x}\,d\nu(x)$
is tight.

We will assume $\nu\in\mathcal{T}_{\epsilon,R}:=\{\nu\in\mathcal{P}_1(M):\nu(\ball{R}{0}^c)>\epsilon\}$ for what follows. We leave $R$ to be defined later, but note that if we define the set $\Omega_R:=\{ \abs{X(0)}>R\}$, then our assumption on $\nu$ yields $\mathbb{P}(\Omega_R)>\epsilon$. By positivity of $L$, this allows us to say that $\mathscr{A}(X)\ge \mathscr{A}(1_{\Omega_R}X)$ (henceforth we define the process $Y(t):=1_{\Omega_R}X(t)$).\\
By (A1), we assume that there is a convex function $\underline{L}:M^\ast\rightarrow\mathbb{R}$ and $C>0$ such that for all $\abs{u}>U$, 
    $    \frac{\underline{L}(u)}{\abs{u}^2}>C.$
Recall that $\underline{L}(\abs{v})$ is a lower bound on $L(t,x,v)$. This imposes a lower bound on the expected action of $Y$:
\begin{equation}\begin{split}\label{eq:Jensencoerc}
    \expect{\int_0^T L(t,Y,\beta_Y(t,Y))\,dt}\ge\,& \expect{\int_0^T \underline{L}(\abs{\beta_Y(t,Y)})\,dt}\overset{(\text{J})}{\ge}\expect{\underline{L}(\abs{V})T}\\
    \overset{(\text{A1})}{>}&C T\expect{1_{\abs{V}>U}\abs{V}^2}\ge C T\bracket{\expect{\abs{V}^2}-U^2},
\end{split}
\end{equation}
where $\beta_Y:=1_{\Omega_R}\beta_X$ is the drift associated with the process $Y$ and $V:=(Y(T)-Y(0))/T$ is its time-average. Hence the expected action of the stochastic process $X$ is bounded:
\begin{equation}\label{eq:Coerctogether}\begin{split}
    \mathscr{A}(X)\ge& \mathscr{A}(Y)> C T\expect{\abs{\frac{Y(T)-Y(0)}{T}}^2 }-CU^2T
    >\frac{C}{T}\expect{\abs{\abs{Y(0)}-\abs{Y(T)}}^2}-CU^2T.
\end{split}
\end{equation}
This leaves us with the same formulation as in (\ref{det.Coercivity.eq}) of the deterministic coercivity result, the remainder of the proof is identical, and the claim is proved. \\
To show that a minimizing sequence $\nu_n$ is sequentially compact in the weak topology, we use  the fact that the set of measures $\nu$ such that
 $   C(\nu,\nu_T)<N\int\abs{y}\,d\nu(y)+\lbstoc(\mu_0,\nu_T)+1$ 
is tight. 
If we let $N:=\int \abs{x}\,d\mu_0(x)$, then the collection of measures such that
\begin{equation*}\begin{split}
    \lbstoc_T^s(\mu_0,\nu_T)+1>&C(\nu,\nu_T)+\underline{W}(\mu_0,\nu)\\
    >&C(\nu,\nu_T)-\int\abs{x}\abs{y}\,d\mu_0(x)\,d\nu(y)\\ 
    \overset{(\text{F})}{=}&C(\nu,\nu_T)-N\int\abs{y}\,d\nu(y),
\end{split}
\end{equation*}
is tight, where (F) is an application of Fubini's theorem. Thus, by Prokhorov's theorem the minimizing sequence of interpolating measures necessarily weakly converges to a minimizing measure.

\begin{rmk}\rm 
a) The same reasoning as in Section 2 yields that 
$C(\nu_0,\nu_1)=C(\nu_1,\nu_0)=\infty$ for $\nu_1\in\mathcal{P}_1(M)$ and $\nu_0\in\mathcal{P}(M)\setminus\mathcal{P}_1(M)$. This implies  that it suffices to take the infimum in (\ref{eq:interpol}) over $\mathcal{P}_1(M)$.
\par
 \noindent b) The attainment of a minimizing interpolating measure $\nu_0$ is sufficient to show the existence of a minimizing $(V,X)$ for $\underline{B}^s_T(\nu_0,\nu_T)$ whenever the latter is finite. This is a consequence of the existence of minimizers for both $\underline{W}(\mu_0,\nu_0)$ and $C^s_T(\nu_0,\nu_T)$ \cite[Proposition 2.1]{M-T}.
\end{rmk}
To establish the duality formula, we will proceed as in the the deterministic case and use the Legendre dual of the optimal cost functional $\nu\to C_T^s(\nu_0, \nu)$, which was derived by Mikami and Thieullen \cite{M-T}. Indeed, they show that 
if the Lagrangian satisfies (A1)-(A3), then 
\begin{equation}\label{Mik}
    C_T^s(\nu_0,\nu_T)=\sup\left\{\int_M f\,d\nu_T-\int_M \Psi^0_{f,-}\,d\nu_0; f\in\mathcal{C}_b^\infty\right\},
\end{equation}
where $\Psi_{f, -}$ is the unique solution to the Hamilton-Jacobi-Bellman equation (\ref{HJB}) that is given by: 
\begin{equation}\label{eq:dynprog}
    \Psi_{_{f,-}}(t,x)=\sup_{X\in\mathcal{A}}\left\{\expcond{f(X(T))-\int_t^T L(s,X(s),\beta_X(s,X))\,ds}{X(t)=x}\right\}.
\end{equation}
Moreover, there exists an optimal process $X$ with drift $\beta_X(t,X)={\rm argmin}_v\{v\cdot \nabla \phi(t,x)+L(t,x,v)\}$.\\
Furthermore, $(\mu,\nu)\mapsto C_T^s(\mu,\nu)$ is convex and lower semi-continuous under the weak$^*$-topology. It follows that 
   $\nu\mapsto\lbstoc_T^s(\mu_0,\nu)$ is weak$^*$-lower semi-continuous on $\mathcal{P}_1(M)$ for all $\mu_0\in\mathcal{P}_1(M^\ast)$, and that 
     $(\mu_0,\nu_T)\mapsto \lbstoc_T^s(\mu_0,\nu_T)$ is jointly convex.
\begin{rmk}\rm 
Note that integrating $\Psi^0_{f,+}$ over $d\nu_0$ yields the Legendre transform of $\nu_T\mapsto C(\nu_0,\nu_T)$ for $f\in\dbound$.
\end{rmk}
For $\mu_0\in \mathcal{P}_1(M^\ast)$, define the function $\underline{B}_{\mu_0}:\mathcal{M}_1(M)\rightarrow\mathbb{R}\cup\{\infty\}$ to be
\begin{equation*}
    \underline{B}_{\mu_0}(\nu):=\begin{cases}\underline{B}(\mu_0,\nu)&\nu\in\mathcal{P}_1(M)\\
    \infty &\text{otherwise}.\end{cases}
\end{equation*}
Since $\underline{B}_{\mu_0}$ is convex and weak$^*$-lower semi-continuous,
we have
\begin{equation}\label{eq:Bddual}
    \underline{B}_{\mu_0}(\nu)=\underline{B}_{\mu_0}^{\ast\ast}(\nu)=\sup_{f\in \lip{M}}\left\{\int f\, d\nu-\underline{B}_{\mu_0}^\ast(f)\right\}.
\end{equation}
We break this into two steps. First we show that when $f\in \dbound$ the dual is appropriate:
\begin{align}\notag
    \underline{B}_{\mu_0}^\ast(f)&\,:=\sup_{\nu_T\in \mathcal{P}_1(M)}\left\{\int f\,d\nu_T-\underline{B}(\mu_0,\nu_T)\right\}\\\label{eq:dualexp}
    &\,\overset{(\ref{eq:interpol})}{=}\sup_{\substack{\nu_T\in\mathcal{P}_1(M)\\\nu\in\mathcal{P}_1(M)}}\left\{\int f\,d\nu_T-C(\nu,\nu_T)-\underline{W}(\mu_0,\nu)\right\}\\\notag
    &\overset{(\ref{eq:dynprog})}{=}\sup_{\nu\in\mathcal{P}_1(M)}\left\{\int \Psi^0_{_{f,-}}(x)\,d\nu(x)-\underline{W}(\mu_0,\nu)\right\}\\\notag
    &\hspace{2mm}= \underline{W}_{\mu_0}^\ast(\Psi^0_{_{f,-}})=-\int (\Psi^0_{f, -})_*\,d\mu_0.
\end{align}
Thus, plugging this into our dual formula (\ref{eq:Bddual}) and restricting our supremum to $\dbound$ gives
\begin{equation*}
    \underline{B}_{\mu_0}(\nu)=\underline{B}_{\mu_0}^{\ast\ast}(\nu)\ge\sup_{f\in \dbound}\left\{\int f\, d\nu+\int (\Psi^0_{f, -})_*\,d\mu_0\right\}.
\end{equation*} 
To show the reverse inequality we will adapt the mollification argument used in \cite[Proof of Theorem 2.1]{M-T}. We assume our mollifier $\eta_\epsilon(x)$ is such that $\eta_1(x)$ is a smooth function on $[-1,1]^d$ that satisfies $\int \eta_1(x)\,dx=1$ and $\int x\eta_1(x)\,dx=0$, then define $\eta_\epsilon(x)=\epsilon^{-d}\eta_1(x/\epsilon)$. Then for Lipschitz $f$, $f_\epsilon:=f\ast \eta_\epsilon$ is smooth with bounded derivatives. We can derive a bound on $\underline{B}_{\mu\ast\eta_\epsilon}^\ast(f)$ by removing the supremum in (\ref{eq:dualexp}) and fixing a process $X\in\mathcal{A}^{\nu_T}$:
\begin{multline*}
    \expect{f_\epsilon(X(T))-\int_0^T L(s,X(s),\beta_X(s,X))\,ds-\langle X(0),V\rangle}\overset{(\text{A2})}{\le}\\
    \expect{f(X(T)+H_\epsilon)-\int_0^T \frac{L(s,X(s)+H_\epsilon,\beta_X(s,X))-\Delta L(0,\epsilon)}{1+\Delta L(0,\epsilon)}\,ds-\langle X(0)+H_\epsilon,V+H_\epsilon\rangle+\abs{H_\epsilon}^2}\le \\
    \frac{D_{\epsilon}^\ast(f\bracket{1+\Delta L(0,\epsilon)})}{1+\Delta L(0,\epsilon)}+T\frac{\Delta L(0,\epsilon)}{1+\Delta L(0,\epsilon)}+d\epsilon^2, 
\end{multline*}
where $D_\epsilon(\nu):=\inf\{(1+\Delta L(0,\epsilon))\underline{W}(\mu_\epsilon,\nu_0)+C(\nu_0,\nu);\nu_0\}$, $H_\epsilon\sim\eta_\epsilon$ is independent of $X(\cdot),V$, thus $X(T)+H_\epsilon\sim d(\eta_\epsilon\ast \nu_T)$. The third line arises by maximizing over processes $(X(\cdot)+H_\epsilon,V+H_\epsilon)$. Note that $\epsilon\mapsto D_\epsilon(\nu)$ is lower semi-continuous for the same reason that $\nu\mapsto \underline{B}_{\mu}(\nu)$ is, and converges to $\underline{B}_{\mu_0}(\nu)$ as $\epsilon\rightarrow 0$. 
\par 
Taking the supremum over $X\in\mathcal{A}_{\mu_0}$ of the left side above, we can retrieve a bound on $\underline{B}_{\mu_0}^\ast(f_\epsilon)$. This bound allows us to say
\begin{equation*}\begin{split}
    \int f_\epsilon\,d\nu-\underline{B}_{\mu}^\ast(f_\epsilon)\ge  \int f\,d\nu_\epsilon-\frac{D_{\epsilon}^\ast(f\bracket{1+\Delta L(0,\epsilon)})}{1+\Delta L(0,\epsilon)}-T\frac{\Delta L(0,\epsilon)}{1+\Delta L(0,\epsilon)}-d\epsilon^2,
\end{split}
\end{equation*}
where we use $\epsilon$-subscript to indicate convolution of a measure with $\eta_\epsilon$. Taking the supremum over $f\in\lip{M}$, we get the reverse inequality:
\begin{equation*}\begin{split}
    \sup_{f\in \dbound}\left\{\int f\,d\nu-\underline{B}_{\mu}^\ast(f)\right\}\ge \frac{D_\epsilon(\nu_\epsilon)}{1+\Delta L(0,\epsilon)}-T\frac{\Delta L(0,\epsilon)}{1+\Delta L(0,\epsilon)}-d\epsilon^2\overset{\epsilon\searrow 0}{\ge}\underline{B}(\mu_0,\nu_T),
\end{split}
\end{equation*}
In the following corollary, we will discuss results pertaining to solutions $\psi_n^t(x):=\psi_n(t,x)$ of the Hamilton-Jacobi-Bellman equation for final conditions $\psi^T_n(x)$. In some sense $\nabla\psi$ is more fundamental than $\psi$, since our dual is invariant under $\psi\mapsto \psi+c$. Thus when discussing the convergence of a sequence of $\psi$, we refer to the convergence of their gradients. Notably the optimal gradient may not be bounded or smooth, hence may not be achieved within the set $\dbound$. In the subsequent corollary, we denote $\mathbb{P}_X$ the measure on $M\times [0,T]$ associated with the process $X$.
\begin{cor} 
\label{cor:minoptX}
Suppose the assumptions on Theorem \ref{thm:stoch.dual}.2 are satisfied and that $\mu_0$ is absolutely continuous with respect to Lebesgue measure. Then $(V,X(t))$ minimizes $\underline{B}(\mu_0,\nu_T)$ if and only if it is a solution to the stochastic differential equation
\begin{align}\label{eq:optX}
    dX =& \nabla_p H(t,X, \nabla \psi(t,X))\, dt + dW_t\\\label{eq:optV} 
    V =& \nabla\bar{\psi}(X(0)),
\end{align}
where $\nabla\psi_n(t,x)\rightarrow \nabla\psi(t,x)$ $\mathbb{P}_X$-a.s. and $\nabla\psi_n(0,x)\rightarrow \nabla\bar{\psi}(x)$ $\nu_0$-a.s. for some sequence $\psi_n(t, x)$ that solves (\ref{eq:HJB}) in such a way that $\psi_n^T:=\psi_n(T, \cdot)$ and $(\psi_n^0)_*:=[\psi_n(0, \cdot)]_*$ are maximixing sequences for the dual problem (\ref{eq:mindual}). Furthermore $\bar{\psi}$ is concave.
\end{cor}
\noindent{\bf Proof:}
First note that there exists such an optimal pair $(V,X)$, in view of Theorem \ref{thm:stoch.dual}.1. Moreover, 
the pair is is optimal iff there exists a sequence of solutions $\psi_n$ to \ref{eq:HJB} that is maximizing in (\ref{eq:mindual})
such that
\begin{equation}\label{eq:converge}
    \expect{\int_0^T L(t,X,\beta_X(t,X))\,dt+\langle X(0),V\rangle}=\lim_{n\rightarrow\infty}\expect{\psi_n^T(X(T))+(\psi_n^0)_*(V)},
\end{equation}
which we can write as
\begin{equation}\label{eq:ineqs}
    \lim_{n\rightarrow\infty}\expect{\underbrace{\psi_n^T(X(T))-\psi_n^0(X(0))}_{\text{(a)}}+\underbrace{\psi_n^0(X(0))-(\psi_n^0)_{**}(X(0))}_{\text{(b)}}+\underbrace{(\psi_n^0)_{**}(X(0))+(\psi_n^0)_*(V)}_{\text{(c)}}},
\end{equation}
where $f_{**}$ is the concave hull of $f$. Applying It\^o's formula to the first two terms, with the knowledge that they satisfy (\ref{eq:HJB}), we get
\begin{equation*}
    \expect{\psi_n^T(X(T))-\psi_n^0(X(0))}=\expect{\int_0^T\langle \beta_X,\nabla\psi_n^t(X(t))\rangle-H(t,X,\nabla\psi_n^t(X(t)))\,dt}.
\end{equation*}
However, by the definition of the Hamiltonian, we have $\langle v,b\rangle - H(t,x,v)\le L(t,x,b)$, which mean that (\ref{eq:ineqs}) yield the following three inequalities:
\begin{align}
   \tag{a} \langle \beta_X,\nabla\psi_n^t(X(t))\rangle-H(t,X,\nabla\psi_n^t(X(t)))\le& L(t,X,\beta_X(t,X))\\
   \tag{b} \psi_n^0(X(0))-(\psi_n^0)_{**}(X(0))\le& 0\\
   \tag{c} (\psi_n^0)_{**}(X(0))+(\psi_n^0)_*(V)\le& \langle V,X(0)\rangle.
\end{align}
In other words, (\ref{eq:ineqs}) breaks the problem into a stochastic and a Wasserstein transport problem (in the flavour of Theorem \ref{thm:stoch.dual}), along with a correction term to account for $\psi^0_n$ not being necessarily concave. Adding (\ref{eq:converge}) to the mix, allows us to obtain $L^1$ convergence in the (a,b,c) inequalities, hence a.s. convergence of a subsequence $\psi_{n_k}$. \\ 
Note that the convergence in (b,c) means that $\psi_n^0$ converges $\nu_0$-a.s. to a concave function $\overline{\psi}$ such that $x\mapsto \nabla\overline{\psi}$ is the optimal transport plan for $\underline{W}(\nu_0,\mu_0)$ \cite{B1}.
\\
To obtain the optimal control for the stochastic process, one needs the uniqueness of the point $p$ achieving equality in (a). This is a consequence of the strict convexity and coercivity of $b\mapsto L(t,x,b)$ for all $t,x$. The differentiability of $L$ further ensures this value is achieved by $p=\nabla_v L(t,x,b)$. Hence (a) holds iff
\begin{equation*}
    \nabla\psi_n^t(X_t)\longrightarrow \nabla_v L(t,X,\beta_X(t,X))\qquad\mathbb{P}_X\text{-a.s.}
\end{equation*}
Since $\psi_n^t$ are deterministic functions, this demonstrates that $X_t$ is a Markov process with drift $\beta_X$ determined by the inverse transform: $\beta_X(t,X)=\nabla_p H(t,X,\nabla\psi(t,X))$, i.e., (\ref{eq:optX}).
\begin{rmk}
It is not possible to conclude from the above work that $\bar{\psi}(x)=\psi(0,x)$ without a regularity result on $t\mapsto\psi(t,x)$ for the optimal $\psi$. This is because $\bar{\psi}$ is defined on a $\mathbb{P}_X$-null set.
\end{rmk}

\section{Deterministic and stochastic Bolza duality}

For the rest of the paper, we shall assume that the Lagrangian $L$ is independent of time, but that it is convex, proper and lower semi-continuous in both variables. We then  consider the dual Lagrangian ${\tilde L}$ defined on $M^*\times M^*$ by 
$$
\tilde L(v,q):=L^*(q, v)=\sup\{\langle v, y\rangle +\langle p, q\rangle -L(y, p);\, (y, p)\in M\times M\}, 
$$
 the corresponding fixed-end costs on $M^*\times M^*$, 
\begin{equation}\label{tilde}
{\tilde c}_T(u, v):=\inf\{\int_0^T{\tilde L}(\gamma (t), {\dot \gamma}(t))\, dt; \gamma\in C^1([0, T), M^*);  \gamma (0)=u, \gamma (T)=v\},
\end{equation}
and its associated optimal transport
\begin{equation}\label{BBT*}
{\tilde C}_T(\mu_0, \mu_T):=\inf \{\int_{M^*\times M^*} {\tilde c}_T(x,y)\, d\pi;\, \pi\in \mK(\mu_0,\mu_T)\}. 
\end{equation}
More specifcally, we shall assume the following conditions on $L$, which are weaker than $(A1), (A2), (A3)$ but for the crucial condition that $L$ is convex in both variables.

\begin{trivlist}

\item (B1) $L:M\times M \to \R\cup\{+\infty\}$ is convex, proper and lower semi-continuous in both variables.

\item (B2) The set $F(x):=\{p; L(x, p) <\infty\}$ is non-empty for all $x\in M$, and for some $\varrho>0$, we have ${\rm dist} (0, F(x)) \leq \varrho (1+|x|)$ for all $x\in M$.

\item (B3) For all $(x, p)\in M\times M$, we have $L(x,p) \geq \theta (\max\{0, |p|-\alpha |x|\})-\beta |x|$, where $\alpha, \beta$ are constants, and $\theta$ is a coercive, proper, non-decreasing function on $[0, \infty)$. 
\end{trivlist}
These conditions on the Lagrangian make sure that the Hamiltonian $H$ is finite, concave in $x$ and convex in $q$, hence locally Lipschitz. Moreover, we have
\begin{equation}
\psi(x)-(\gamma |x| +\delta)|q| \leq H(x, q) \leq \phi(q)+(\alpha |q| +\beta)|x| \hbox{for all $x, q$ in $M\times M^*$,}
\end{equation}
where $\alpha, \beta, \gamma, \delta$ are constants, $\phi$ is finite and convex and $\psi$ is finite and concave (see \cite{R-W2}.\\
We note that under these conditions, the cost $(x, y)\to c_t(x,y)$ is convex proper and lower semi-continuous on $M\times M$. But the cost $b_T$ is nicer in many ways. For one, it is everywhere finite and locally Lipschitz continuous on $[0, \infty)\times M\times M^*$. However, the main addition in the case of joint convexity for $L$ is the following so-called Bolza duality that we briefly describe in the deterministic case since it had been studied in-depth in various articles by T. Rockafellar \cite{R1} and co-authors \cite{R-W1, R-W2}. The stochastic counterpart is more recent and has been established by Boroushaki and Ghoussoub \cite{B-G}.\\ 
We consider the path space 
${\mathcal A}^{2}_{M}: ={\mathcal A}^{2}_{M}[0, T]= \{ u:[0,T] \rightarrow M; \,\dot{u} \in L^{2}_{M}  \}
$
equipped with the norm
\[
     \|u\|_{{\mathcal A}^{^{2}}_{M}} = \left(\|u(0)\|_{M}^{2} +
     \int_0^T \|\dot{u}\|^{2} dt\right)^{\frac{1}{2}}.
     \]
     Let $L$ be a convex Lagrangian on $M\times M$ as above, $\ell$ be a proper convex lower semi-continuous function on $M\times M$ and consider the minimization problems,
\begin{equation}\label{Bolza}
\hbox{$({\mathcal P})\qquad \inf\left\{
\int_0^TL(\gamma (s), {\dot \gamma}(s))\, ds +\ell (\gamma(0), \gamma (T)); \, \gamma \in  C^1([0, T), M)\right\},$}
\end{equation}
and  
\begin{equation}
\hbox{$({\tilde {\mathcal P}})\qquad \inf\left\{\int_0^T{\tilde L}(\gamma (s), {\dot \gamma}(s))\, ds +\ell^* (\gamma(0), -\gamma (T)); \, \gamma \in  C^1([0, T), M)\right\}$.}
\end{equation}
\begin{thm} Assume $L$ satisfies (B1), (B2) and (B3), and that $\ell$ is proper, lsc and convex.
\begin{enumerate}
\item  If there exists $\xi$ such that $\ell (\cdot, \xi)$ is finite, or there exists $\xi'$ such that $\ell (\xi',\cdot)$ is finite, then $$\inf({\mathcal P})= -\inf({\tilde {\mathcal P}}).$$
 This value is not $+\infty$, and if it is also not  $-\infty$, then there is an optimal arc $v(t)\in {{\mathcal A}}^2[0, T]$ for $({\tilde {\mathcal P}})$. 
\item A similar statement holds if we replace $\ell$ by $\tilde \ell$ in the above hypothesis and $({\tilde{\mathcal P}})$ by $({\mathcal P})$ in the conclusion. 
\item If both conditions are satisfied, then both $({\tilde{\mathcal P}})$ and $({\mathcal P})$ are attained respectively by optimal arcs $v(t), x(t)$ in ${{\mathcal A}}^2[0, T]$. 
\end{enumerate}
\end{thm}
In this case, these arcs satisfy
$(\dot v(t), v(t))\in \partial L(x(t), \dot x(t))$ for a.e. $t,$
which can also be written in a dual form 
$({\dot x}(t), x(t))\in \partial {\tilde L}(v(t), \dot v(t))$ for a.e. $t$,
or in a Hamiltonian form as
\begin{eqnarray}
\dot x(t)&\in& \partial_v H(x(t), v(t))\\
-\dot v(t)&\in& {\tilde \partial}_x H(x(t), v(t)), 
\end{eqnarray} 
coupled with the boundary conditions
\begin{equation}
(v(0), -v(T))\in \partial \ell (x(0), x(T)).
\end{equation}
See for example \cite{R1}. The above duality has several consequences. 
\begin{prop} The value function $\Phi_{g, +}(x)=\inf\{g(y)+c_t(y, x);\, y\in M\}$, which is the variational solution of the Hamilton-Jacobi equation (\ref{HJ+}) starting at $g$, can be expressed in terms of the  $b$ and $\tilde c$ costs as follows:

\begin{enumerate}
\item If $g$ is convex and lower semi-continuous, then $\Phi_{g, +}(t,x)=\sup\{b_t( v, x)-g^*(v); \, v\in M^*\}$ is convex lower semi-continuous for every $t\in [0, +\infty)$.

\item The convex Legendre transform of $\Phi_{g, +}$ is given by the formula
$$
\tilde\Phi_{g^*, +}(t, w)=\inf\{g^*(v)+{\tilde c}_t(v, w);\, v\in M^*\}.
$$

\item For each $t$, the graph of the subgradient $\partial \Phi_{g, +}(t, \cdot)$, i..e., 
$\Gamma_g(t)=\{(x, v); v\in \partial \Phi_{g, +}(t, x)\}$
is a globally Lipschitz manifold of dimension $n$ in $M \times M^*$, which depends continuously on $t$.

\item If a Hamiltonian trajectory $(x(t), v(t))$ over $[0, T ]$ starts with $v(0) \in \partial g(x(0))$, then $v(t) \in \partial \Phi_{g, +} (t, x(t))$ for all $t \in [0, T]$. Moreover, this happens 
 if and only if
$x(t)$ is optimal in the minimization problem  that defines $\Phi_{g, +}(t, x)$ and
$v(t)$ is optimal in the minimization problem that defines $\tilde\Phi_{g^*, +}(t, w)$.
\end{enumerate}
\end{prop}
\begin{rmk}\rm
The above shows that in the case when $L$ is jointly convex, the corresponding forward Hamilton-Jacobi equation has convex solutions whenever the initial state is convex, while the corresponding backward Hamilton-Jacobi equation has concave solutions if the final state is concave. Unfortunately, we shall see that in the mass transport problems we are considering, one mostly propagates concave (resp., concave) functions forward (resp., backward), hence losing their concavity (resp., convexity). 

This said, the cost functionals  $c_T$, ${\tilde c}_T$, $b_T$ are all value functions $\Phi_g$ starting or ending with affine function $g$. Indeed, $b_t( v, x)=\Phi_{g, +}(t,x)$, when $g_v(y)=\langle v, y\rangle$. In this case, $g_v^*(u)=0$ if $u=v$  and $+\infty$ if $u\neq v$, which yields that the Legendre dual of $x\to \Phi_{g, +}(t, x)=b_t( v, x)$ is $w\to {\tilde c}_t(v , w)$. One can also deduce the following.
\end{rmk}

\begin{prop} Under assumptions $(B_1), (B_2), (B_3)$ on the Lagrangian $L$, the costs $c$ and $b$ have the following properties:
\begin{enumerate}

\item For each $t\geq 0$, $(x, y)\to c_t(x,y)$ is convex proper and lower semi-continuous on $M\times M$.

\item For each $t\geq 0$, $v\to b_t( v, x)$ is concave on $M^*$, while $x\to b_t( v, x)$ is convex  on $M$. Moreover, $b$ is locally Lipschitz continuous on $[0, \infty)\times M\times M^*$.

\item  The costs $b$, $c$ and $\tilde c$ are dual to each other in the following sense:
\begin{itemize}
\item For any $(v,x)\in M^*\times M$, we have $b_t(v, x)=\inf\{\langle v, y\rangle +c_t(y,x);\,  y\in M\}.$

\item For any $(y, x)\in M\times M$, we have $c_t(y,x)=\sup\{ b_t( v, x)-\langle v, y\rangle; v\in M^*\}.$

\item For any $(v,x)\in M^*\times M$, we have $b_t(v, x)=\sup\{\langle w, x\rangle - {\tilde c}_t(v, w); w\in M^*\}. $
\end{itemize}
\item The following properties are equivalent:
\begin{enumerate}
\item $(-v,w)\in \partial_{y,x}c_T(y, x)$;

\item $w \in \partial_{x}b_T (v, x)$ and $y \in {\tilde \partial}_{v}b_T(v, x)$.

\item There is a Hamiltonian trajectory $(\gamma (t), \eta (t))$ over $[0, T]$ starting at $(y, v)$ and ending at $(x, w)$.

\end{enumerate}
 
\end{enumerate}
\end{prop}
     This leads us to the following standard condition in optimal transport theory.

\begin{defn} A cost function $c$ satisfies {\it the twist condition} if for each $y\in M$, we have $x=x'$ whenever the differentials  $\partial_yc(y,x)$ and $\partial_yc(y,x')$ exist and are equal.
\end{defn}
In view of the above proposition, $c_T$ satisfies the twist condition if there is at most one Hamiltonian trajectory starting at a given initial state $(v, y)$, while the cost $b_T$ satisfies the twist condition if for any given states $(v, w)$, there is at most one Hamiltonian trajectory starting at $v$ and ending at $w$.  

\subsection*{The stochastic Bolza duality and its applications}

We now deal with the stochastic case. 
We define the {\it It\^o space} $\mathcal{I}^p_M$ consisting of all $M$-valued processes of the following form:
\begin{align}\label{spaceA22}
\begin{split}
\mathcal{I}^p_M= \Big\{X: &\Omega_T \rightarrow M ; \  X(t) = X_0 + \int_0^t \beta_X(s) ds+ \int_0^t \sigma_X(s) dW(s), \\
&\textnormal{for} \ X_0 \in  L^2(\Omega,\F_0,\P;M), \ \beta_X \in L^p(\Omega_T; M), \, \sigma_X \in L^2(\Omega_T; M)\Big\},
\end{split}
\end{align}
where $\beta_X$ and $\sigma_X$ are both progressively measurable and $\Omega_T:=\Omega\times[0,T]$. The cases of $p=1,2,\infty$ will be of interest to us. We equip $\mathcal{I}^2_M$ with the norm
\begin{align*}
\Vert X \Vert^2_{\mathcal{I}^2_M}= \E \left(\Vert X(0)\Vert^2_M +\int_0^T \Vert \beta_X(t) \Vert^2_{M} \, dt+ \int_0^T  \Vert \sigma_X(t)\Vert^2_{M} \, dt  \right), 
\end{align*}
so that it becomes a Hilbert space. 
 The dual space $(\mathcal{I}^2_M)^*$ can also be identified with $L^2(\Omega; M) \times L^2(\Omega_T;M) \times  L^2(\Omega_T; M)$. In other words, each $q \in (\mathcal{I}^2_M)^*$ can be represented by the triplet 
$$q = (q_0, q_1(t), Q(t)) \in L^2(\Omega; M) \times L^2(\Omega_T; M) \times  L^2(\Omega_T; M),$$ 
in such a way that the duality can be written as: 
\begin{equation}\label{duality_A22}
\langle X, q \rangle_{{\mathcal{I}^2_M}\times (\mathcal{I}^2_M)^*} =\E   \Big\{ \langle q_0,X(0) \rangle_M + \int_0^T\langle q_1(t),\beta_X(t) \rangle_M\, dt + \frac{1}{2} \int_0^T\langle Q(t),\sigma_X(t) \rangle_{M} \, dt \Big \}.
\end{equation}
Similarly, the dual of $\mathcal{I}^1_M$ can be identified with $\mathcal{I}^\infty_M$.\par

We shall use the following result recently established in \cite{B-G}.
\begin{thm}\label{partial_sd_thm}{\bf (Boroushaki-Ghoussoub)}
Let $(\Omega, \F, \F_t, P)$ be a  complete probability space with normal filtration, and let $L(\cdot, \cdot)$ and $M$ be two jointly convex Lagrangians on $M \times M$, 
Assume $\ell$ is a convex lsc function on $M \times M$. Consider the Lagrangian on $\mathcal{I}^2_M\times (\mathcal{I}^2_M)^*$ defined by
\begin{align}\label{bolza.lagrangian}
\begin{split}
\mathcal{L}(X,p)&=
\E \, \Big\{\int_0^T L(X(t)-p_1 (t), -\beta_X(t)) \, dt+ \ell(X(0)-p_0,X(T)) \\ &\quad + \ \frac{1}{2} \int_0^T M(\sigma_X(t) -P(t),-\sigma_X(t)) \, dt \Big\}. 
\end{split}
\end{align}
Its Legendre dual is then given for each $q:=(0,q_1,Q)$ by 
\begin{align*}
\mathcal{L}^*(q,Y)&=   \E \Big\{\ell^*(-Y(0),Y(T)) +   \int_0^T  L^*(-\beta_Y(t), Y(t)-q_1(t)) \, dt \\
&\quad + \frac{1}{2}\,  \int_0^T  M^*(-\sigma_Y(t),\sigma_Y(t)-Q(t)) \, dt\Big\}.
\end{align*}
\end{thm}
Note that standard duality theory implies that in general
\begin{equation}
    \inf_{X\in \mathcal{I}^2}\{\mathcal{L}(X,0)\}\ge \sup_{Y\in (\mathcal{I}^2)^\ast}\{-\mathcal{L}(0,Y)\}.
\end{equation}
In our case we shall restrict ourselves to processes of fixed diffusion. This facilitates the proving of a stochastic analog to Bolza duality:
\begin{prop}
Assume $L$ satisfies (A1) and (A2), and there exists (a.s.-)unique $V_0\in L^2(\mathbb{P})$ such that $\ell^\ast(V_0,\cdot)<\infty$ and (a.s.)-unique $\sigma_V\in L^2(\mathbb{P}\times\lambda_{[0,T]})$ such that $M^\ast(\sigma_V,\cdot)<\infty$, then there is no duality gap, ie. 
\begin{equation}
    \inf_{X\in \mathcal{I}^2}\{\mathcal{L}(X,0)\}= \sup_{Y\in (\mathcal{I}^2)^\ast}\{-\mathcal{L}^\ast(0,Y)\}
\end{equation}
\end{prop}
Note that, unlike the deterministic case, there there is no backwards condition that works if there is an $V_T\in L^2(\mathbb{P})$ such that $\ell^\ast(\cdot,V_T)<\infty$, this is because stochastic processes, in general, are irreversible.\\
{\bf Proof:}
We begin with augmenting our space by considering $\beta_V\in L^1(\mathbb{P}\times\lambda_{[0,T]})$---we call this augmented set $\mathcal{I}^1$. If we can show the duality gap is satisfied in $\mathcal{I}^1$, by our coercivity condition (A2) we can then show that it must be satisfied in $\mathcal{I}^2$.\\
 We proceed by a variational method outlined by Rockafellar \cite{R1}.
First, we define
\begin{equation}
    \phi(q):=\inf_{Y\in (\mathcal{I}^1)^\ast}\{\mathcal{L}^\ast(q,Y)\}.
\end{equation}
As the infimum of a jointly convex function, $\phi$ itself is convex. The benefit of this definition is that
\begin{equation}
    \phi^\ast(X)=\sup_{q,v}\{\langle X,q\rangle-\mathcal{L}^\ast(q,v) \}=\mathcal{L}^{\ast\ast}(X,0)=\mathcal{L}(X,0).
\end{equation}
Hence,  $X$ minimizes $\mathcal{L}$ if and only if
\begin{equation}
    X\in\partial\phi(0)\quad\Longleftrightarrow\quad \phi(0)+\phi^\ast(X)=0\quad\Longleftrightarrow\quad \mathcal{L}(X,0)=-\inf_{Y\in(\mathcal{I}^2)^\ast}\{{\cal L}^* (0, Y)\}.
\end{equation}
In other words, there is no duality gap if and only if $\partial \phi(0)$ is non-empty. Note that this holds if there is an open (relative to $\{q; \phi(q)<\infty\}$) neighbourhood $N$ of the origin in $\mathcal{I}^\infty$ such that $\mathcal{L}^\ast(q,Y)<\infty$ for $q\in N$.\\ 
By our assumptions, we may fix $Y_0, \sigma_Y$ to be the unique elements such that $\ell(Y_0,\cdot)<\infty$ and $M^\ast(\sigma_Y,\cdot)<\infty$ (guaranteeing subdifferentiability in these variables), and let $Y=(Y_0,\beta_Y,\sigma_Y)$ be such that $\mathcal{L}^\ast(0,Y)<\infty$. For a perturbation $\beta_V\in L^\infty(\mathbb{P}\times\lambda_{[0,T]})$ with $\norm{\beta_V}_\infty<\epsilon$, note that (A2) gives for all $(t,u)\in [0,T]\times M^\ast$, 
\begin{equation}
    L(t,Y_t-\beta_V,u)<\bracket{1+\Delta L(0,\epsilon)}L(t,Y_t,u)+\Delta L(0,\epsilon),
\end{equation}
and
\begin{equation}\begin{split}
    \phi(V)=&\inf_{Y\in\mathcal{I}^1_M}\mathcal{L}^\ast(V,Y)\\
    \le &\E\ell^\ast(-Y_0,Y_T)+\E\int_0^T \tilde{L}(t,Y_t-\beta_V(t),\beta_Y(t))\,dt\\
    \le &\E\ell^\ast(-Y_0,Y_T)+(1+\Delta\tilde{L}(0,\epsilon))\E\int_0^T \tilde{L}(t,Y_t,\beta_Y(t))\,dt+T\Delta\tilde{L}(0,\epsilon),
\end{split}
\end{equation}
which is finite for $\norm{\beta_V}_\infty<\epsilon$ sufficiently small by (A2). Hence $\phi$ is finite and continuous in a open set of the origin (all relative to its domain), and duality is achieved on $\mathcal{I}^1$.\par 
To show that this duality is achieved in $\mathcal{I}^2$, it suffices to remark that $\E\int_0^T L(t,Y_t,\beta_Y)\ge \E \int\underline{L}(\beta_Y)\,dt\ge C\E \int\abs{\beta_Y}^2-B\,dt=\infty$ for $\beta_Y\in L^1(\mathbb{P}\times\lambda_{[0,T]})\setminus L^2(\mathbb{P}\times\lambda_{[0,T]})$ (where $C,B$ are fixed constants).

\section{Maximizing the ballistic cost: Deterministic case}

With Bolza duality in mind, we can now turn to the maximizing ballistic cost.

\begin{thm} \label{interpol} Assume that $L$ satisfies hypothesis (B1), (B2) and (B3), and let $\nu_T$ be a probability measure with compact support on $M$, that is also absolutely continuous with respect to Lebesgue measure. Then,  
\begin{enumerate}
\item The following interpolation formula holds:
\begin{equation}\label{HL.two}
{\overline B}_T(\mu_0,\nu_T)=\sup\{ {\overline W}(\nu_T, \mu)-{\tilde C}_T(\mu_0, \mu);\, \mu\in {\mathcal P}(M^*)\}.  
\end{equation}
The supremum is attained at some probability measure $\mu_T$ on $M^*$, and the final Kantorovich potential for ${\tilde C}_T(\mu_0, \mu_T)$ is convex. 

 \item We also have the following duality formulae: 
\begin{equation}\label{ballistic.dual2}
{\overline B}_T(\mu_0,\nu_T)=\inf\left\{\int_M h(x)\, d\nu_T(x)+ \int_{M^*} \tilde\Phi^0_{h^*,-}(v)\, d\mu_0(v); \,  \hbox{$h$ convex in $\lip{M}$} \right\}. 
\end{equation}
 and
 \begin{equation}\label{ballistic.dual3}
{\overline B}_T(\mu_0,\nu_T)=\inf\left\{\int_M ({\tilde\Phi}^T_{g,+})^*(x)\, d\nu_T(x)+ \int_{M^*} g(v)\, d\mu_0(v); \,  \hbox{$g$ in Lip(M$^*$)} \right\}. 
\end{equation}

\item There exists a convex function $h: M^* \to \R$ such that 
\begin{equation}
{\overline B}_T(\mu_0,\nu_T)= \int _{M^*} b_T \large({S}^*_T \circ \nabla h^*(x), x\large) d\nu_T(x),
\end{equation}
where 
${S}^*_T(v)=\pi^*\phi^{H_*}_T(v, \nabla  h),$
and $\phi^{H_*}_t$ the flow  associated to the Hamiltonian $H_*(v,x)=-H(-x, v)$, whose Lagrangian is $L_*(v,q)=L^*(-q, v)$. 
In other words, an optimal map for ${\overline B}_T(\mu_0,\nu_T)$ is given by the inverse of the map $x\to \pi^*\phi^{H_*}_T(\nabla h^* (x), x)$.

\item We also have 
\begin{equation}
{\overline B}_T(\mu_0,\nu_T)= \int _{M^*} b_T (v, \nabla h \circ {\tilde S}_Tv\large) d\mu_0(v),
\end{equation}
where 
$
{\tilde S}_T(v)=\pi^*\phi^{{\tilde H}}_T(v, {\tilde d}_vh_0),
$ and 
$\phi^{{\tilde H}}_t$ being the Hamiltonian flow associated to $\tilde L$ (i.e., ${\tilde H}(v, x)=-H(x, v)$, and $h_0=\tilde\Phi^0_{h^*,-}$.

$h_0$ the solution $h(0, v)$ of the backward Hamilton-Jacobi equation (\ref{dHJ}) with $h(T, v)=h(v)$. 
\end{enumerate}
\end{thm}
{\bf Proof:} To show (\ref{HL.two}) and (\ref{ballistic.dual2}), first note that for any probability measure $\mu$ on $M^*$, we have 
\begin{equation}\label{inequality}
{\overline B}_T(\mu_0,\nu_T)\geq  {\overline W}(\nu_T, \mu)-{\tilde C}_T(\mu_0, \mu).  
\end{equation}
Indeed, since $\nu_T$ is assumed to be absolutely continuous with respect to Lebesgue measure, Brenier's theorem yields a convex function $h$ that is differentiable $\mu_T$-almost everywhere on $M$ such that $(\nabla h)_\#\nu_T=\mu$, and  ${\overline W}(\nu_T, \mu)=\int_{M}\langle x, \nabla h(x)\rangle \, d\nu_T(x)$. Let $\pi_0$ be an optimal transport plan for 
${\tilde C}_T(\mu_0, \mu)$, that is $\pi_0\in {\mathcal K}(\mu_0, \mu)$ such that
$
{\tilde C}_T(\mu_0, \mu)=\int_{M^*\times M^*} {\tilde c}_T(v,w)\, d\pi_0(v, w).
$
Let ${\tilde \pi}_0:= S_\#\pi_0$, where $S(v,w)=(v, \nabla h^*(w))$, which is a transport plan in ${\mathcal K}(\mu_0, \nu_T)$. Since $ b_T(v, y)\geq \langle \nabla h(x), y\rangle -{\tilde c}_T(v, \nabla h(x) )$ for every $(y, x, v)\in M\times M\times M^*$, we have 
\begin{eqnarray*}
{\overline B}_T(\mu_0,\nu_T)&\geq&\int_{M^*\times M} b_T(v,x)\, d{\tilde \pi}_0(v, x)\\
&\geq& \int_{M^*\times M}\{ \langle \nabla h(x), x\rangle -{\tilde c}_T(v, \nabla h(x) )\} d{\tilde \pi}_0(v, x)\\
&= &\int_{M}\langle x, \nabla h(x)\rangle \, d\nu_T(x) -\int_{M^*\times M^*}{\tilde c}_T(v, w)\} d\pi_0(v, w)\\
&=& {\overline W}(\nu_T, \mu)-{\tilde C}_T(\mu_0, \mu).  
\end{eqnarray*}
To prove the reverse inequality, we use standard Monge-Kantorovich theory to write
\begin{eqnarray*}
{\overline B}_T(\mu_0,\nu_T)&=&\sup\big\{\int_{M^*\times M} b_T(v, x)\, d\pi(v, x);\, \pi \in {\mathcal K}(\mu_0, \nu_T)\big\}\\
&=&\inf\big\{\int_{M} h(x)\, d\nu_T(x)-\int_{M^*}g(v)\, d\mu_0(v);\, h(x)-g(v) \geq b_T(v, x)\big\},
\end{eqnarray*} 
where the infimum is taken over all admissible Kantorovich pairs $(g, h)$ of functions, i.e. those satisfying
the relations
$$
g(v)=\inf_{x\in M}h(x)- b_T(v, x)
\text{\quad  and \quad } 
h(x)=\sup _{v\in M^*} b_T(v, x))+g(v)
$$
Note that $h$ is convex.
Since the cost function $b_T$ is continuous, the supremum ${\overline B}_T(\mu_0,\nu_T)$ is attained at some probability measure $\pi_0\in {\mathcal K}(\mu_0, \nu_T)$. Moreover, the infimum in the dual problem  is attained at some pair $(g, h)$ of admissible Kantorovich functions. It follows that $\pi_0$ is supported on the set
\[
{\mathcal O}:=\{(v, x)\in M^*\times M; \, b_T(v, x)=h(x)-g(v)\}.
\]
We now exploit the convexity of $h$, and use the fact that 
for each $(v, x) \in {\mathcal O}$, the function $y\to h(y)-g(v)-b_T(v,y)$ attains its minimum at $x$, which means that 
$
\nabla h(x)\in {\partial}_x b_T(v, x).
$
But since $\tilde c_T$ is the Legendre transform of $b_T$ with respect to the $x$-variable, we then have 
\begin{equation}
b_T(v, x)+{\tilde c}_T(v, \nabla h(x))=\langle x, \nabla h(x)\rangle \,\, \hbox{on ${\mathcal O}$}.
\end{equation}
Integrating with $\pi_0$, we get since $\pi_0\in {\mathcal K}(\mu_0, \nu_T)$,
\begin{equation}
\int_{M^*\times M}b_T(v, x)\, d\pi_0+\int_{M^*\times M}{\tilde c}_T(v, \nabla h(x))d\pi_0=\int_{M}\langle x, \nabla h(x)\rangle\, d\nu_T.
\end{equation}
Letting $\mu_T=\nabla h_{\#}\nu_T$, we obtain that 
\begin{equation}\label{crucial}
{\overline B}_T(\mu_0,\nu_T)+\int_{M^*\times M}{\tilde c}_T(v, \nabla h(x))d\pi_0={\overline W}(\nu_T, \mu_T), \end{equation}
where 
$
{\overline W}(\nu_T, \mu_T)=\sup\{\int_{M\times M^*} \langle x,v\rangle\, d\pi;\, \pi \in {\mathcal K}(\nu_T, \mu_T)\}.
$
Note that we have used here that $h$ is convex to deduce that ${\overline W}(\nu_T, \mu_T)=\int_{M}\langle x, \nabla h(x)\, d\mu_T$ by the uniqueness in Brenier's decomposition. We now prove that 
\begin{equation}\label{magic}
\int_{M^*\times M}{\tilde c}_T(v, \nabla h(x))d\pi_0={\tilde C}_T(\mu_0, \mu_T). 
\end{equation}
Indeed, we have 
$\int_{M^*\times M}{\tilde c}_T(v, \nabla h(x))d\pi_0\geq {\tilde C}_T(\mu_0, \mu_T)$
since the measure $\pi=S_\#\pi_0$, where $S(v,x)=(v, \nabla h(x))$ has marginals $\mu_0$ and $\mu_T$ respectively. On the other hand, 
 (\ref{crucial}) yields 
\begin{eqnarray*}
\int_{M^*\times M}{\tilde c}_T(v, \nabla h(x))d\pi_0&=&\int_{M}\langle x, \nabla h(x)\rangle\, d\nu_T(x)-\int_{M^*\times M}b_T(v, x)\, d\pi_0\\
&=&\int_{M}h^*(\nabla h(x)) d\nu_T(x)+ \int_Mh(x)\, d\nu_T(x)+\int_{M^*}g(v)\, d\mu_0(v)-\int_Mh(x)\, d\nu_T(x)\\
&=&\int_{M^*}h^*(w) d\mu_T(w)+\int_{M^*}g(v)\, d\mu_0(v).
\end{eqnarray*}
Moreover, since $h(x)-g(v) \geq b(v, x)$, we have  
$
h^*(w) +g(v) \leq {\tilde c}_T(v, w).  
$
Indeed, since for any $(v, w)\in M^*\times M^*$, we have 
${\tilde c}(t,v,w)=\sup\{\langle w, x\rangle - b_t( v, x); x\in M\}, $
it follows that for any $y\in M$, 
\begin{eqnarray*}
{\tilde c}_T(v,w)\geq  \langle w, y\rangle - b_t( v, y)\geq \langle w, y\rangle +g(v)-h(y), 
\end{eqnarray*}
hence
$h^*(w) +g(v) \leq {\tilde c}_T(v, w)$, 
which means that the couple $(-g, h^*)$ is an admissible Kantorovich pair for the cost ${\tilde c}_T$. Hence,
\begin{eqnarray*}
{\tilde C}_T(\mu_0, \mu_T)&\leq& \int_{M^*\times M}{\tilde c}_T(v, \nabla h(x))d\pi_0\\
&=&\int_{M}h^*(w) d\mu_T(w)+\int_{M^*}g(v)\, d\mu_0(v)\\
&\leq&\sup\{\int_{M^*} \phi_T(w)\, d\mu_T(w)-\int_{M^*}\phi_0(v)\, d\mu_0(v) ;\,  \phi_T(w)-\phi_0(v)  \leq {\tilde c}_T(v, w)\}\\
&=&{\tilde C}_T(\mu_0, \mu_T).
\end{eqnarray*}
It follows that
${\overline B}_T(\mu_0,\nu_T)={\overline W}(\nu_T, \mu_T)-{\tilde C}_T(\mu_0, \mu_T)$. 
In other words, the supremum in (\ref{inequality}) is attained by the measure $\mu_T$. Note that the final optimal Kantorovich potential for ${\tilde C}_T(\mu_0, \mu_T)$ is $h^*$, hence is convex.\\
The first duality formula (\ref{ballistic.dual3}) follows since we have established that if $(g, h)$ are an optimal pair of Kantorovich functions for ${\overline B}_T(\mu_0,\nu_T)$,  then $(g, h^*)$ are an 
optimal pair of Kantorovich functions for ${\tilde C}_T(\mu_0, \mu_T)$. In other words,  the initial Kantorovich function for ${\overline B}_T(\mu_0,\nu_T)$ is $g=\tilde\Phi_{h^*, -}(0, \cdot)$. This proves formula (\ref{ballistic.dual2}).\\
 To show  (\ref{ballistic.dual3}), we can --now that the interpolation (\ref{HL.two}) is established--proceed as in Section 2, by identifying the Legendre transform of the functionals $\nu\to {\overline W}(\nu, \nu_T)$ and 
$\mu \to {\tilde C}_T(\mu, \mu_T)$.

To show part 3), we start with the interpolation inequality and write that 
\[
{\overline B}_T(\mu_0,\nu_T)={\overline W}(\nu_T, \mu_T)-{\tilde C}_T(\mu_0, \mu_T),
\]
for some probability measure $\mu_T$. The proof also shows that there exists a convex function $h: M^\ast \to \R$ and another function $k:M^\ast\to \R$ such that $(\nabla h)_\#\mu_T=\nu_T$,
$
{\overline W}(\nu_T, \mu_T)=\int_{M} \langle \nabla h(v), v\rangle d\mu_T(v).
$
and
$
\tilde{C}_T(\mu_0, \mu_T)=\int_{M^\ast}h(u)\, d\mu_T(u)-\int_{M^\ast} k(v)\, d\mu_0(v).
$
Now use the theorem of Fathi-Figalli to write 
\begin{equation}
\tilde{C}_T(\mu_0, \mu_T)=\int _{M^\ast} c_T (v, \tilde{S}_Tv) d\mu_0(v), \end{equation}
where $\tilde{S}_T(v)=\pi^*\phi^{\tilde{H}}_T(v,\tilde d_v k)$. 
Note that 
\begin{equation}
\hbox{${\overline B}_T(\mu_0,\nu_T)\geq \int _{M^*} b_T (v, \nabla h\circ \tilde{S}_T (v)) d\mu_0(v),$}
\end{equation}
since $(\tilde{S}_T)_\#\mu_0=\mu_T$ and $\nabla h_\#\mu_T=\nu_T$, and  therefore $(I \times \nabla h\circ \tilde{S}_T)_\#\mu_0$ belongs to ${\mathcal K}(\mu_0, \nu_T)$.\\
On the other hand, since $ b_T(u, x)\geq  \langle \nabla h(v), x\rangle-\tilde{c}_T(u,\nabla h(v))$ for every $v\in M^*$, we have
\begin{eqnarray*}
{\overline B}_T(\mu_0,\nu_T)&\geq& \int _{M^*} b_T (v, \nabla h\circ \tilde{S}_T (v)) d\mu_0(v)\\
&\geq& \int_{M^*}\{\langle \nabla h\circ \tilde{S}_T(v), \tilde{S}_T(v)\rangle-\tilde{c}_T(v,\tilde{S}_T(v)) \}\, d \mu_0(v)\\
&=&\int _{M^\ast} \langle\nabla h(v), v\rangle d\mu_T(v) - \int_{M^*}\tilde{c}_T(v,\tilde{S}_T(v)) \, d \mu_0(v)   \\
&=& {\overline W}(\nu_T, \mu_T)-\tilde{C}_T(\mu_0, \mu_T)\\
&=&{\overline B}_T(\mu_0,\nu_T).
\end{eqnarray*}
It follows that 
${\overline B}_T(\mu_0,\nu_T)=\int _{M^*} b_T (v, \nabla h\circ \tilde{S}_T (v)) d\mu_0(v).$\\
To get (3), use the pushforward $\nu_T=(\nabla h\circ \tilde{S}_T )_\#\mu_0$ to write the above in terms of the measure $\nu_T$, using the fact that $(\nabla h)^{-1}=\nabla {h^\ast}$ and $\tilde{S}_T^{-1}=S^\ast_T$ where $S^\ast_T(v)=\pi^\ast\phi^{H_\ast}_t(v,\tilde{d}_vh)$ and $\phi^{H_\ast}_t$ is the Hamiltonian flow associated to the hamiltonian $H_\ast(v,x):=-H(-x,v)$. This gives us
\begin{equation*}
{\overline B}_T(\mu_0,\nu_T)=\int _{M^*} b_T (S_T^\ast\circ\nabla {h^\ast}(x), x) d\nu_T(x)=\int _{M^*} b_T (\pi^\ast\phi^{H_\ast}_t(\nabla {h^\ast}(x),\tilde{d}_vh), x) d\nu_T(x). 
\end{equation*}
Since $h$ is convex, we have that $\tilde d_xh=\nabla h(x)$, hence $\tilde d_{\nabla {h^\ast}(x)}h=\nabla h\circ \nabla {h^\ast}(x)=x$, which yields our claim that 
$$
 {\overline B}_T(\mu_0,\nu_T)=\int _{M} b_T \large( \pi^*\phi^{\tilde{H}}_T(\nabla {h^\ast}(x), x)\large,x) d\nu_T(x).
 $$

\begin{rmk}\rm
While the costs $c$ and $\tilde{c}_T$ are themselves jointly convex in both variables, one cannot deduce much in terms of the convexity or concavity of the corresponding Kantorovich potentials. 
However, we note that the interpolation (\ref{HL.one}) of  ${\underline B}_T(\mu_0,\nu_T)$ selects a $\nu_0$ such that ${C}_T(\nu_0, \nu_T)$ has a concave initial Kantorovich potential, while the interpolation (\ref{HL.two}) of  ${\overline B}_T(\mu_0,\nu_T)$  
selects a $\mu_T$ such that ${\tilde C}_T(\mu_0, \mu_T)$ has a convex final Kantorovich potential.

Furthermore, one wonders whether the formula 
\begin{equation}
c_t(y,x)=\sup\{ b_t( v, x)-\langle v, y\rangle; v\in M^*\},
\end{equation}
also extends to Wasserstein space. We show it under the condition that the initial Kantorovich potential of $C_T(\nu_0, \nu_T)$ is concave, and conjecture that it is also a necessary condition.
\end{rmk}

\begin{thm}\label{endpoints} Assume $M=\R^d$ and that $L$ satisfies hypothesis (B1), (B2) and (B3). Assume $\nu_0$ and $\nu_T$ are probability measures on $M$ such that $\nu_0$ is absolutely continuous with respect to Lebesgue measure. If the initial Kantorovich potential of $C_T(\nu_0, \nu_T)$ is concave then the following holds: 
\begin{equation}\label{three}
C_T(\nu_0, \nu_T)=\sup\{{\underline B}_T(\mu,\nu_T)-{\underline W}(\nu_0, \mu);\, \mu\in {\mathcal P}(M^*)\},  
\end{equation}
and the supremum is attained.
\end{thm}
\noindent {\bf Proof:} Again, it is easy to show that 
\begin{equation}\label{three.prime}
C_T(\nu_0, \nu_T)\geq \sup\{{\underline B}_T(\mu,\nu_T)- {\underline W}(\nu_0, \mu);\, \mu\in {\mathcal P}(M^*)\}.  
\end{equation}
To prove equality,  
we assume  
that the initial Kantorovich potential $g$ is concave and write 
\begin{eqnarray*}
C_T(\nu_0,\nu_T)&=&\inf\{\int_{M\times M} c(y, x)\, d\pi(y, x);\, \pi \in {\mathcal K}(\nu_0, \nu_T)\}\\
&=&\sup\{\int_{M} h(x)\, d\nu_T(x)-\int_{M}g(y)\, d\nu_0(y);\, h(x)-g(y) \leq c_T(y, x)\}.
\end{eqnarray*} 
Since the cost function $c_T$ is continuous, the infimum $C_T(\nu_0,\nu_T)$ is attained at some probability measure $\pi_0\in {\mathcal K}(\nu_0, \nu_T)$. Moreover, the infimum in the dual problem  is attained at some pair $(g, h)$ of admissible Kantorovich functions.  It follows that $\pi_0$ is supported on the set
\[
{\mathcal O}:=\{(y, x)\in M \times M; \, c_T(y, x)=h(x)-g(y)\}
\]
Since $g$ is concave, use the fact that  for each $(y, x) \in {\mathcal O}$, the function $z\to h(x)-g(z) -c_T(z,x)$ attains its maxmum at $y$, to deduce  that 
$
-\nabla g(y)\in {\partial}_y c_T(y, x).$\\
Since $g$ concave and  $b_t(v, x)=\inf\{\langle v, z\rangle +c_t(z,x);\,  z\in M\}$, this means that for $(y, x) \in {\mathcal O}$,
\begin{equation}
c_T(y, x)=b_T(\nabla g(y), x)-\langle \nabla g(y), y\rangle.
\end{equation}
Integrating with $\pi_0$, we get since $\pi_0\in {\mathcal K}(\nu_0, \nu_T)$,
\begin{equation}
\int_{M\times M}c_T(y, x)\, d\pi_0=\int_{M\times M}b_T(\nabla g(y), x)\, d\pi_0-\int_{M}\langle \nabla g(y), y\rangle\, d\nu_0.
\end{equation}
Letting $\mu_0=(\nabla g)_{\#}\nu_0$, and since $g$ is concave, we obtain that 
\begin{equation}\label{crucial.2}
C_T(\nu_0,\nu_T)=\int_{M\times M}b_T(\nabla g(y), x)\, d\pi_0 -{\underline W}(\nu_0, \mu_0).
 \end{equation}
We now prove that 
\begin{equation}\label{magic.2}
\int_{M\times M}b_T(\nabla g(y), x)\, d\pi_0(y,x)={\underline B}_T(\mu_0, \nu_T). 
\end{equation}
Indeed, we have 
$\int_{M\times M}b_T(\nabla g(y), x)\, d\pi_0\geq {\underline B}_T(\mu_0, \nu_T),$
since the measure $\pi=S_\#\pi_0$ where $S(y, x)=(\nabla g(y), x)$ has $\mu_0$ and $\nu_T$ as marginals.
On the other hand, 
(\ref{crucial.2}) yields 
\begin{eqnarray*}
\int_{M\times M}b_T(\nabla g(y), x)\, d\pi_0&=&\int_{M\times M}c_T(y, x)\, d\pi_0+\int_{M}\langle y, \nabla g(y)\rangle\, d\nu_0(y)\\
&=&\int_Mh(x)\, d\nu_T(x)-\int_{M}g(y)\, d\nu_0(y)-\int_{M}(-g)^*(-\nabla g(y)) d\nu_0(y)+ \int_Mg(y)\, d\nu_0(y)\\
&=&\int_{M}h(x)\, d\nu_T(x)-\int_{M^*}(-g)^*(-v) d\mu_0(v).
\end{eqnarray*}
Moreover, since $h(x)-g(y) \leq c_T(y, x)$, it is easy to see that 
$h(x)-(-g)^*(-v) \leq b_T(v,x)$, that is the couple $((-g)^*(-v), h (x))$ is an admissible Kantorovich pair for the cost $b_T$. It follows that
\begin{eqnarray*}
{\underline B}_T(\mu_0, \nu_T)&\leq&\int_{M\times M}b_T(\nabla g(y), x)\, d\pi_0\\
&=&\int_{M}h(x)\, d\nu_T(x)-\int_{M}(-g)^*(-v) d\mu_0(v)\\
&\leq&\sup\{\int_{M} \phi_T(x)\, d\mu_T(x)-\int_{M^*}\phi_0(v)\, d\mu_0(v);\,  \phi_T(x)-\phi_0(v)  \leq b_T(v, x)\}\\
&=&{\underline B}_T(\mu_0, \nu_T),
\end{eqnarray*}
and 
$C_T(\nu_0, \nu_T)={\underline B}_T(\mu_0,\nu_T)-{\underline W}(\nu_0,  \mu_0). $
In other words, the supremum in (\ref{three}) is attained by the measure $\mu_0$.

\begin{cor} Assume $M=\R^d$ and that $L$ satisfies hypothesis (B1), (B2) and (B3). Assume $\nu_0$ and $\nu_T$ are probability measures on $M$ such that $\nu_0$ is absolutely continuous with respect to Lebesgue measure, and that the initial Kantorovich potential of $C_T(\nu_0, \nu_T)$ is concave. If $b_T$ satisfies the twist condition, then there exists a map $X_0^T: M^*\to M$ and a concave function $g$ on $M$ such that 
\begin{equation}
C_T(\nu_0, \nu_T)=\int _{M} c_T (y, X_0^T\circ \nabla g(y)) d\nu_0(y).
\end{equation}
\end{cor}
\noindent{\bf Proof:} In this case,
$C_T(\nu_0, \nu_T)={\
\underline B}_T(\mu_0,\nu_T)-{\underline W}(\nu_0, \mu_0), $
for some probability measure $\mu_0$ on $M^*$. Let $g$ be the concave function on $M$ such that 
$(\nabla g)_\#\nu_0=\mu_0$ and 
$
{\underline W}(\nu_0, \mu_0)=\int_M \langle \nabla g(y), y\rangle d\nu_0(y).
$
Since $b_T$ satisfies the twist condition, there exists a map $ X_0^T: M^*\to M$  such that $(X_0^T)_\#\mu_0=\nu_T$ and 
\begin{equation}
{\underline B}_T(\mu_0, \nu_T)=\int _{M^*} b_T (v, X_0^Tv) d\mu_0(v).
\end{equation}
Note that the infimum $C_T(\nu_0, \nu_T)$  is attained at some probability measure $\pi_0\in {\mathcal K}(\nu_0, \nu_T)$ and that $\pi_0$ is supported on a subset ${\mathcal O}$ of $M\times M$ such that 
for $(y, x) \in {\mathcal O}$,
$c_T(y, x)=b_T(\nabla g(y), x)-\langle \nabla g(y), y\rangle.$
Moreover, $C_T(\nu_0,\nu_T)=\int_{M\times M}b_T(\nabla g(y), x)\, d\pi_0 -{\underline W}(\nu_0, \mu_0)$, and 
\[
\int_{M\times M}b_T(\nabla g(y), x)\, d\pi_0={\underline B}_T(\mu_0, \nu_T)=\int _{M^*} b_T (v, X_0^Tv) d\mu_0(v)=\int _{M} b_T (\nabla g(y), X_0^T\circ \nabla g(y)) d\nu_0(y).
\]
Since $b_T$ satisfies the twist condition, it follows that for any $(y, x)\in {\mathcal O}$, we have that 
$x=X_0^T\circ \nabla g(y)$ from which follows that  $C_T(\nu_0, \nu_T)=\int _{M} c_T (y, X_0^T\circ \nabla g(y)) d\nu_0(y).$

\begin{cor} \label{cor1}Consider the cost $c_1(y, x)=c(x-y)$, where $c$ is a convex function on $M$ and let $\nu_0, \nu_1$ be  probability measures on $M$ such that the initial Kantorovich potential associated to $C_T(\nu_0, \nu_T)$ is concave. Then, there exist concave functions $\phi: M\to \R$, $\psi: M^*\to \R$ and a probability measure $\mu_0$ on $M^*$ such that
\begin{equation}
(\nabla \psi\circ \nabla \phi)_\#\nu_0=\nu_1,
\end{equation}
and 
\begin{equation}
C_1(\nu_0, \nu_1)+\int_{M^*}c^*(v)\, d\mu_0(v)=\int _{M} c ( \nabla \psi\circ \nabla \phi(y)-y) d\nu_0(y)=\int_M \langle \nabla {\psi_*}(y)-\nabla \phi (y), y\rangle \, d\nu_0(y). 
\end{equation}
\end{cor}
\noindent {\bf Proof:} The cost $c(x-y)$ corresponds to $c_1(y, x)$, where the Lagrangian is $L(x, v)=c(v)$, that is 
\begin{equation}\label{v}
c_1(y,x)=\inf\{\int_0^1c({\dot \gamma}(t))\, dt; \gamma\in C^1([0, 1), M);  \gamma(0)=y, \gamma(1)=x\}=c(x-y).
\end{equation}
It follows from (\ref{three}) that there is a probability measure $\mu_0$ on $M^*$ such that 
$C_1(\nu_0, \nu_1)={\underline B}_1(\mu_0,\nu_1)-{\underline W}(\nu_0, \mu_0)$. But
in this case,
$b_1(v, x)=\inf\{\langle v, y\rangle +c(x-y);\,  y\in M\}=\langle v, x\rangle-c^*(v)$, hence 
\begin{equation}\label{ng}
C_1(\nu_0, \nu_1)={\underline B}_1(\mu_0,\nu_1)-{\underline W}(\nu_0, \mu_0)={\underline W}(\mu_0, \nu_1)-\int_{M^*}c^*(v)\, d\mu_0(v)-{\underline W}(\nu_0, \mu_0). 
\end{equation}
In other words, 
\begin{equation}\label{ng}
C_1(\nu_0, \nu_1)+K={\underline W}_1(\mu_0,\nu_1)-{\underline W}(\nu_0, \mu_0),
\end{equation}
where $K$ is the constant $\int_{M^*}c^*(v)\, d\mu_0(v)$. 

Apply Brenier's theorem twice to find concave functions $\phi: M\to \R$ and $\psi: M^*\to \R$ such that $(\nabla \phi)_\#\nu_0=\mu_0$, $(\nabla \psi)_\#\mu_0=\nu_1$ and
\[
\hbox{${\underline W}(\nu_0, \mu_0)=\int_{M}\langle y, \nabla \phi(y)\rangle \, d \nu_0(y)$
\quad and \quad 
$
{\underline W}( \mu_0, \nu_1)=\int_{M^*}\langle v, \nabla \psi(v)\rangle \, d\mu_0(v).$}
\]
It follows from the preceeding corollary that 
\begin{eqnarray*}
C_1(\nu_0, \nu_1)+K=\int _{M} c_1 (y, \nabla \psi\circ \nabla \phi(y)) d\nu_0(y)=\int _{M} c ( \nabla \psi\circ \nabla \phi (y)-y) d\nu_0(y).
\end{eqnarray*}
Note also that 
\begin{eqnarray*}
C_1(\nu_0, \nu_1)+K&=&\int_{M}\langle v, \nabla \psi(v)\, d \mu_0(v) -\int_{M}\langle y, \nabla \phi(y)\, d \nu_0(y)
\\
&=&\int_M \langle \nabla {\psi_*}(y), y\rangle \, d\nu_0(y)-\int_{M}\langle y, \nabla \phi(y)\, d \nu_0(y)\\
&=& \int_M \langle \nabla {\psi_*}(y)-\nabla \phi (y), y\rangle \, d\nu_0(y). 
\end{eqnarray*}

\section{Maximizing the ballistic cost: Stochastic case}

Define the transportation cost between two random variables $V$ on $M^*$ and $X$ on $M$ by:
\begin{equation}
    b^s_T(V,Y):=\inf\{\expect{\langle V,X(0)\rangle +\int_0^T L(t,X_t,\beta_X(t))\,dt}; X\in\mathcal{A}, X(T)=Y\text{ a.s.}\},
\end{equation}
where $\mathcal{A}$ indicates It\^o processes with Brownian diffusion.  The minimizing ballistic cost considered earlier is then
\begin{equation}
    \underline{B}^s_T(\mu_0,\nu_T)=\inf\{b^s_T(V,Y); V\sim \mu_0, Y\sim\nu_T\}, 
\end{equation}
while the maximizing cost is defined as:
\begin{equation}\label{stoch.sup}
    \overline{B}^s_T(\mu_0,\nu_T):=\sup\{b^s_T(V,Y); V\sim \mu_0, Y\sim\nu_T\}.
\end{equation}
\begin{thm}\label{max.Stoch.dual} Assume $L$ is a Lagrangian on $M\times M^*$ such that $L$ and its dual $\tilde{L}$ satisfies (A0)-(A3).
\begin{enumerate}
\item  The following formula then holds:
\begin{equation}\label{eq:Bover}
    \overline{B}_T^s(\mu_0,\nu_T):=\sup\left\{\expect{\langle X,V(T)\rangle-\int_0^T {\tilde L}(t,V,\beta_V(t,V))\,dt}; V\in \mathcal{A}, V_0\sim\mu_0, X\sim \nu_T\right\}.
\end{equation}

\item The following duality holds:
\begin{equation}\label{eq:maxinterpol}
    \overline{B}_T^s(\mu_0,\nu_T)=\sup\{\overline{W}(\mu,\nu_T)-\tilde{C}^s_T(\mu_0,\mu);\, \mu \in\mathcal{P}_1(M^\ast)\},
\end{equation}
where $\tilde{C}^s_T$ is the action corresponding to the Lagrangian $\tilde L$. Furthermore, if $\nu_0\in\mathcal{P}_1(M)$, and $\mu_T\in\mathcal{P}_1(M^\ast)$ there exist an optimal interpolant $\mu_T$ in $\mathcal{P}_1(M^\ast).$ 

\item If $\mu_0\in\mathcal{P}_1(M^*)$, $\nu_T$ has compact support, and $\overline{B}(\mu_0,\nu_T)<\infty$, then 
\begin{equation}\label{eq:maxdual}
    \overline{B}_s(\mu_0,\nu_T)=\inf\left\{\int_{M} g\,d\nu_T+\int_{M^*}\tilde \Psi^0_{g^*, -}\,d\mu_0;\,  \hbox{$g \in \dbound (M)$ and convex} \right\},
\end{equation}
where $\tilde \Psi_{g^*, -}$ solves the Hamilton-Jacobi-Bellman equation on $M^*$
\begin{align}\label{eq:HJB2}\tag{HJB2}
    \pderiv{\psi}{t}+\frac{1}{2}\Delta \psi -H(t,\nabla_v\psi, v)=0 \qquad \psi(v,T)=g^*(v).
\end{align}
\end{enumerate}
\end{thm}
\noindent{\bf Proof:} 1) For a fixed pair $(V, Y)$, we consider the Bolza energy $\mathcal{L}_{(V,Y)}$ --defined in (\ref{bolza.lagrangian})-- associated to $L$ and the two Lagrangians $\ell$ and $M$ defined as:
\begin{align}
    \ell_{(Y,U)}(\omega,y,z):=&\begin{cases}\langle z,U(\omega)\rangle &y=Y(\omega)\\
    \infty &\text{else}
    \end{cases}&
    M(\xi,\zeta):=&\begin{cases}-\zeta-1 &\xi=1\\
    \infty &\text{else}.
    \end{cases}
\end{align}
Note that the minimizing stochastic cost can be written as,
\begin{equation}
    \underline{B}_s(\mu_0,\nu_T):=\inf\{\inf\{\mathcal{L}_{(V,Y)}(X(t),0);X\in \mathcal{I}^2\};V\sim\mu_0,Y\sim \nu_T\}
\end{equation}
while the maximizing cost is 
\begin{equation}
    \overline{B}_s(\mu_0,\nu_T)=\sup\{\inf\{\mathcal{L}_{(V,Y)}(X(t),0);X\in \mathcal{I}^2\};V\sim\mu_0,Y\sim\nu_T\}.
\end{equation}
Applying Bolza duality turns the infimum to a supremum:
\begin{equation}
    \overline{B}_s(\mu_0,\nu_T)=\sup\{\sup\{-\mathcal{L}_{(V,Y)}^\ast(0,U(t)); U\in \mathcal{I}^2\};V\sim\mu_0,Y\sim\nu_T\},
\end{equation}
which results in (\ref{eq:Bover}).\\
\noindent 2) The proof of the interpolation result can now follow closely  the proof for the minimization problem.\\
3) We again try to identify the Legendre transforms of the functionals $\nu\mapsto \overline{W}(\mu,\nu)$ and $\mu \to \tilde{C}^s_T(\mu_0,\mu).$ We obtain easily that 
\begin{itemize}
\item If $\mu\in\mathcal{P}_1(M^*)$ has compact support, then for all $f\in {\rm Lip}(M)$, then 
\begin{equation*}
    \sup_{\nu\in\mathcal{P}_1(M)}\left\{\int_{M} f\,d\nu+\overline{W}(\mu,\nu)\right\}=\int_{M^*} (-f)^\ast\,d\mu.
\end{equation*}

\item If  $g \in\dbound(M^*)$, then 
\[
\sup_{\mu\in\mathcal{P}_1(M^\ast)} \left\{\int_{M^*} g\,d\mu -\tilde{C}^s_T(\mu_0,\mu)\right\}=\int_{M^*}\tilde\Psi_{g,-}\, d\mu_0.
\]
\end{itemize}
Define $\overline{B}_{\mu_0}:\nu\mapsto\overline{B}(\mu_0,\nu)$, and note that the interpolation formula (\ref{eq:maxinterpol}) and a result of Mikima-Thieullen \cite{M-T} concerning $\tilde C_T^s$ yields that $\overline{B}_{\mu_0}$ is a concave function. Furthermore it is weak$^\star$-upper semi-continuous on $\mathcal{P}_1(M)$. Thus we have
\begin{equation}\label{eq:maxddual}
    \overline{B}_{\mu_0}(\nu_T)=-(-\overline{B}_{\mu_0})^{\ast\ast}(\nu_T)=\inf_{f\in \lip{M}}\left\{-\int_{M} f\,d\nu_T+(-\overline{B}_{\mu_0})^\ast(f)\right\}.
\end{equation}
Investigating the dual, we find
\begin{align}\notag
    (-\overline{B}_{\mu_0})^\ast(f)=&\sup_{\nu\in\mathcal{P}_1(M)} \left\{\int_{M} f\,d\nu +\overline{B}_{\mu_0}(\nu)\right\}\\\notag
    =&\sup_{\substack{\mu\in\mathcal{P}_1(M^\ast)\\\nu\in\mathcal{P}_1(M)}} \left\{\int_{M} f\,d\nu +\overline{W}(\mu, \nu)-\tilde{C}^s_T(\mu_0,\mu)\right\}\\ \label{eq:lipmaxdual}
    =&\sup_{\mu\in\mathcal{P}_1(M^\ast)} \left\{\int_{M^*} (-f)^\ast\,d\mu -\tilde{C}^s_T(\mu_0,\mu)\right\}.
\end{align}
Note that in the case where $(-f)^\ast\in\dbound$, this is simply $\int_{M^*}\tilde\Psi_{(-f)^*,-}\, d\mu_0$, yielding
\begin{equation*}
    \overline{B}_{\nu_0}(\mu_T)\le \inf_{(-f)^\ast\in \dbound}\left\{-\int_{M^*} f\,d\mu_T+\int_{M^*}\tilde\Psi_{(-f)^*,-}\, d\mu_0\right\}.
\end{equation*}
In either case, we can restrict our $f$ to be concave by noting that if we fix $g=(-f)^\ast$, then the set of corresponding $\{-f; (-f)^\ast=g\}$ is minimized by the convex function $g^\ast=(-f)^{\ast\ast}\le -f$ \cite[Proposition 4.1]{Ekeland}. Thus it suffices to consider $f$ convex.\par 
We now show that it is sufficient to consider this infimum over convex $g\in\dbound$ by a similar mollification argument to that used for $\underline{B}$ (note that the mollifying preserves convexity). Maintaining the same assumptions and notation as in our earlier argument, we first note a useful application of Jensen's inequality to the legendre dual of a mollified function:
\begin{equation*}\begin{split}
    g_\epsilon^\ast(v)=&\sup_x\left\{\langle v,x\rangle -\expect{g(x+H_\epsilon)}\right\}
    \overset{(\text{J})}{\le}\sup_x\left\{\langle v,x\rangle -g(x)\right\}= g^\ast(v).
\end{split}
\end{equation*}
Mikami \cite[Proof of Theorem 2.1]{M-T} further shows that 
\begin{equation*}
    (\ref{eq:lipmaxdual})=C^\ast_{\nu_0}(g_\epsilon)\le \frac{C^\ast_{\nu_0\ast\eta_\epsilon}((1+\Delta L(0,\epsilon))g)}{1+\Delta L(0,\epsilon)}+T\frac{\Delta L(0,\epsilon)}{1+\Delta L(0,\epsilon)}.
\end{equation*}
Putting these together we get
\begin{equation*}
    \int g^\ast_\epsilon\,d\mu_T+(-\overline{B}_{\nu_0})^\ast(g_\epsilon^\ast)\,d\nu_0\le \int g^\ast\,d\mu_T+\frac{C^\ast_{\nu_0\ast\eta_\epsilon}((1+\Delta L(0,\epsilon))g)}{1+\Delta L(0,\epsilon)}+T\frac{\Delta L(0,\epsilon)}{1+\Delta L(0,\epsilon)}.
\end{equation*}
And once we take the infimum over convex $g\in\lip{M}$, we get
\begin{equation*}
    \inf\left\{\int g^\ast\,d\mu_T+\bracket{-\overline{B}_{\nu}}^\ast(-g^\ast);g\text{ convex in }\dbound\right\}\le \frac{-(-\overline{B})^{\ast\ast}_{\nu_0\ast\eta_\epsilon}(\mu_{L,\epsilon})}{1+\Delta L(0,\epsilon)}+T\frac{\Delta L(0,\epsilon)}{1+\Delta L(0,\epsilon)},
\end{equation*}
where $d\mu_{L,\epsilon}(v):=d\mu_T(\bracket{1+\Delta L(0,\epsilon)}v)$. Taking $\epsilon\searrow 0$ dominates the right side by $\overline{B}(\mu_0,\nu_T)$ (where we exploit the upper semi-continuity of $\overline{B}$), completing the reverse inequality.
\begin{cor}[Optimal Processes for $\overline{B}$]\label{cor:maxoptX}
Suppose the assumptions on Theorem \ref{max.Stoch.dual} are satisfied, with $\mu_0$ absolutely continuous with respect to Lebesgue measure. Then, the pair $(V,X)$ is optimal for (\ref{stoch.sup}) if and only there is an Ito process $V(t)$ that satisfy the backward Stochastic differential equation,
\begin{align}\label{eq:maxoptX}
    dV =& \nabla_p H(t,\nabla \psi(t,V),V)\, dt + dW_t\\\label{eq:maxoptV} 
    X =& \nabla\bar{\psi}(V(T)),
\end{align}
where $\lim_{n\rightarrow \infty}\psi_n(T,x)\rightarrow \bar{\psi}(x)$ $\nu_T$-a.s. and $\lim_{n\rightarrow\infty}\psi_n(t,x)=\psi(t,x)$ $\mathbb{P}_V$-a.s. for some sequence $\psi_n(t, x)$ that solves (\ref{eq:HJB}) in such a way that $\psi_n^0=\psi_n(0, \cdot)$ and $\psi_n^T=\psi_n(T, \cdot)$ are a minimizing pair for the dual problem.
\end{cor}
\noindent{\bf Proof:}
If $(V,X)$ is optimal, then Theorem \ref{max.Stoch.dual} means there exists a sequence of solutions $\psi_n(t, v)$ to (\ref{eq:HJB}) with convex final condition $\psi_n^T$, such that
\begin{eqnarray}\label{eq:maxconverge}
    \expect{\langle X,V(T)\rangle-\int_0^T {\tilde L}(t,V,\beta_V(t,V))\,dt}&=&\lim_{n\rightarrow\infty}\expect{\bracket{\psi_n^T}^\ast(X)+\psi_n^0(V(0))},
\end{eqnarray}
which we write as
\begin{equation*}
    \lim_{n\rightarrow\infty}\expect{\bracket{\psi_n^T}^\ast(X)+\psi_n^T(V(T))-\psi_n^T(V(T))+\psi_n^0(V(0))}.
\end{equation*}
Applying It\^o's formula to the last two terms, with the knowledge that $\psi_n$ satisfies (\ref{eq:HJB}) we get
\begin{equation*}
    \expect{-\psi_n^T(V(T))+\psi_n^0(V(0))}=\expect{\int_0^T-\langle \beta_V,\nabla\psi_n^t(V(t))\rangle-H(t,\nabla\psi_n^t(V(t)),V(t))\,dt}
\end{equation*}
However, by the definition of the Hamiltonian, we have $-\langle q,v\rangle - H(t,x,v)\ge -\tilde{L}(t,v,q)$, similarly $\psi^\ast(v)+\psi(x)\ge\langle v,x\rangle$. These inequalities allow us to separate the limit in (\ref{eq:maxconverge})  
into two requirements:\\
 (a) $\langle \beta_V,\nabla\psi_n^t(V(t))\rangle+H(t,\nabla\psi_n^t(V(t)),V(t))$ must converge to $\tilde{L}(t,V,\beta_V(t,V))$ and \\
 (b) $\psi_n^T(V(T))+\bracket{\psi_n^T}^\ast(X)$ must converge to $\langle X,V(T)\rangle$ in $L^1$ hence a subsequence $\psi_{n_k}$ exists such that this convergence is a.e.\par  
The journey from (a) to (\ref{eq:maxoptX}) 
is as in Corollary \ref{cor:minoptX}. 
The only difference from the earlier corollary is that we know that $\psi_n$ must converge to a convex function, so (b) implies $X=\nabla\lim_{n\rightarrow \infty}\psi_n(V(T))$.

 \section{Final Remarks}
 
 The interpolation formula can be seen as a Hopf-Lax formula on Wasserstein space, since for a fixed $\mu_0$ on $M^*$ (resp., fixed  $\nu_T$ on $M$), then as a function of the terminal (resp., initial) measure, we have 
 \begin{equation}
{\underline {\mathcal B}}^{\mu_0}(t, \nu)=\inf\{{\underline {\mathcal U}}^{\mu_0}(\varrho)+ C_t(\varrho, \nu);\, \varrho \in {\mathcal P}(M)\}\,\, \hbox{and \,\, ${\overline {\mathcal B}}^{\nu_T}(t, \mu)=\inf\{{\overline {\mathcal U}}^{\nu_T}(\varrho)- {\tilde C}_t(\varrho, \mu);\, \varrho\in {\mathcal P}(M^*)\},$ }
\end{equation}
where 
$${\underline {\mathcal U}}^{\mu_0}(\varrho)={\underline W}(\mu_0, \varrho)\hbox{\quad and \quad ${\overline {\mathcal U}}^{\nu_T}(\varrho)={\overline W}(\nu_T, \varrho).$}
$$  
 The following Eulerian formulation illustrates best how ${\underline {\mathcal B}}^{\mu_0}(t, \nu)$ and ${\overline {\mathcal B}}^{\nu_T}(t, \mu)$
 can be represented as value functionals on Wasserstein space. Indeed,  lift the Lagrangian $L$ to the tangent bundle of Wasserstein space via the formula 
\[
\hbox{${\mathcal L}(\varrho, w);=\int_M L(x, w(x)) \, d\varrho(x)$\quad and \quad $\tilde {\mathcal L}(\varrho, w);=\int_{M^*} {\tilde L}(x, w(x)) \, d\varrho(x),$
}
\]
where $\varrho$ is any probability density  on $M$ (resp., $M^*$) and $w$ is a vector field on $M$ (resp., $M^*$). 
\begin{cor} Assume $L$ satisfies hypothesis (A0) and (A1), 
and let $\mu_0$  be a probability measure on $M^*$ with compact support, then  
\begin{eqnarray}
{\underline {\mathcal B}}^{\mu_0}(T, \nu):={\underline B}_T(\mu_0,\nu)
&=& \inf\left\{{\underline {\mathcal U}}^{\mu_0}(\varrho_0) +
\int_0^T {\mathcal L}(\varrho_t, w_t) dt;\,  \partial_t \varrho+ \nabla \cdot (\varrho w)=0,\,  \varrho_T=\nu\right\},  
 \end{eqnarray}
 The set of pairs $(\varrho, w)$ considered above are such that $t \rightarrow \varrho_t \in \mathcal P_1(M)$ (resp.,  $t \rightarrow w_t (x)\in {\rm Lip}(\R^n))$ are paths of Borel fields.

\end{cor}
One can then ask whether these value functionals also satisfy a Hamilton-Jacobi equation on Wasserstein space such as
  \begin{equation}\label{e:gangbo:master}
\left\{
\begin{array}{ll}
& \partial_t B 
+ {\mathcal H}(t, \nu, \nabla_\nu B(t, \nu)) 
=0,\\ [5pt] 
& B(0, \nu) ={\underline W}(\mu_0, \nu). 
\end{array}
\right.
\end{equation}  
Here the Hamiltonian  is defined as 
\[
{\mathcal H}(\nu, \zeta) =\sup\{\int \langle \zeta, \xi\rangle d\nu -{\mathcal L}(\nu, \xi); \xi \in T_\nu^*({\mathcal P}(M))\}.
\]
We note that Ambrosio-Feng \cite{A-F} have shown recently that --at least in the case where the Hamiltonian is the square-- value functionals on Wasserstein space yield a unique {\it metric viscosity solution} for (\ref{e:gangbo:master}). As importantly, Gangbo-Sweich \cite{G-S} have shown recently that under certain conditions, value functionals yield solutions to the so-called {\it Master equations} of mean field games. We refer to their paper for the relevant definitions.

\begin{thm} (Gangbo-Swiech) Assume $\mathcal U_0: \mathcal P(M) \rightarrow \mathbb R$, and $U_0:  M \times \mathcal P(M) \rightarrow \mathbb R$ are functionals such that 
$ \nabla_x U_0(x, \mu) \equiv \nabla_\mu \mathcal U_0(\mu)(x)$ for all $x \in M$, $\mu \in \mathcal P(M),$ and consider the value functional,
\begin{equation*} 
\mathcal U(t, \nu)=\inf \left\{ \mathcal U_0(\varrho_0) + \int_0^t \mathcal L(\varrho, w) dt;\,  \partial_t \varrho+ \nabla \cdot (\varrho w)=0,\,  \varrho_T=\nu\right\}.
 \end{equation*} 
Then, there exists $U: [0,T] \times M \times \mathcal P(M) \rightarrow \mathbb R$ such that 
\begin{equation*}\label{1b}
 \nabla_x U_t(x, \nu) \equiv \nabla_\nu \mathcal U_t(\nu)(x) \quad \hbox{ for all $x \in M$, $\nu \in \mathcal P(M),$}
\end{equation*}  
 and  $U$ satisfies the Master equation below (\ref{e:gangbo:master.1}).
\end{thm} 
Applied to the value functional ${\underline {\mathcal B}}^{\mu_0}(t, \nu):={\underline B}_t(\mu_0,\nu)$, this 
should then yield the existence, for any probabilities $\mu_0, \nu_T$, of a function $\beta: [0,T] \times M \times \mathcal P(M) \rightarrow \mathbb R$ such that 
\begin{equation*}
 \nabla_x \beta (t, x, \nu) \equiv \nabla_\nu {\underline {\mathcal B}}^{\mu_0}(t, \nu)(x) \quad \hbox{ for all $x \in M$, $\nu \in \mathcal P(M),$}
\end{equation*}  
and  $\varrho \in AC^2((0, T)\times {\mathcal P}(M))$ such that 
\begin{equation}\label{e:gangbo:master.1}
\left\{
\begin{array}{ll}
& \partial_t \beta  + \int\langle \nabla_\nu \beta (t, x, \nu)\cdot \nabla  H(x, \nabla_x\beta)\rangle\, d\nu
+ H(x, \nabla_x \beta(t, x,\nu)) =0,\\ [5pt] 
&\partial_t\varrho +\nabla (\varrho \nabla H(x, \nabla_x \beta))=0, \\ [5pt]
& \beta(0, \cdot, \cdot) = \beta_0, \quad \varrho (T, \cdot)=\nu_T,
\end{array}
\right.
\end{equation}  
where $\beta_0(x, \varrho)=\phi_\varrho(x)$, where $\phi_\varrho$ is the convex function such that $\nabla \phi_\varrho$ pushes $\mu_0$ into $\varrho$.

We may furthermore derive a Eulerian formulation of the minimizing stochastic problem. Recall that in Corollary \ref{cor:minoptX} we showed that the optimal process for the minimizing stochastic cost is Markovian. Hence its drift may be described by a vector field, allowing an Eulerian formulation of the process:
\begin{cor}
Assume $L$ satisfies the assumptions (A0)-(A3), then
\begin{eqnarray}
    \underline{\mathcal{B}}^{\mu_0}(T,\nu):&=&\underline{B}^s_T(\mu_0,\nu)\nonumber\\&=&\inf\{\underline{\mathcal{U}}^{\mu_0}(\varrho_0)+\int_0^T\mathcal{L}(t,\varrho_t,w_t)\,dt;\, \, \partial_t\varrho+\nabla\cdot(w\varrho)+\frac{1}{2}\Delta\varrho=0,\varrho_T=\nu\}.
\end{eqnarray}
\end{cor}
\noindent{\bf Proof:}
It can be seen by It\^o's formula that $(\varrho,w)$ is a solution (in the sense of distributions) to $\partial_t\varrho+\nabla\cdot(w\varrho)+\frac{1}{2}\Delta\varrho=0$ iff $\varrho_t\in\mathcal{P}(M)$ is the law of $X_t$ where $X$ solves
\begin{equation}
    X_t=X_0+\int_0^t w(s,X_s)\,ds+W_t.
\end{equation}
Hence the above Eulerian formulation is equivalent to the stochastic process formulation in the case where the optimal drift is described by a Borel vector field. Corollary \ref{cor:minoptX} shows this is the case for $\underline{B}_s$. \hfill $\Box$\\

Finally, we mention that one would like to consider value functionals on Wasserstein space that are more general than those starting with the Wasserstein distance. One can still obtain such functionals via mass transport by considering more general ballistic costs of the form 
\begin{equation}
b_g(T, v, x):=\inf\left\{g(v, \gamma (0)) +\int_0^TL(\gamma (t), {\dot \gamma}(t))\, dt; \gamma \in C^1([0, T), M)\right\},    
\end{equation}
where $g: M^*\times M \to \R$ is a suitable function.

\end{document}